\documentclass[11pt]{amsart}
\usepackage[margin=2.7cm]{geometry}
\usepackage{amsmath,amssymb,amsthm,hyperref,xcolor,tikz}
\newtheorem{theorem}{Theorem}[section]
\newtheorem{lemma}[theorem]{Lemma}
\newtheorem{proposition}[theorem]{Proposition}
\newtheorem{corollary}[theorem]{Corollary}
\newtheorem{conjecture}[theorem]{Conjecture}
\newtheorem{problem}[theorem]{Problem}
\theoremstyle{definition}
\newtheorem{definition}[theorem]{Definition}
\theoremstyle{remark}
\newtheorem{remark}[theorem]{Remark}

\newcommand{\PD}{\mathrm{PD}}
\newcommand{\PDO}{\mathrm{PDO}}

\title{Partitions with designated summands: $k$-tuples and their congruences}

\author{Daniel Yaqubi$^*$}
\address{Department of Computer Science, University of Torbat-e Jam, Torbat-e Jam, Iran}
\email{yaqubi@tjamcaas.ac.ir, daniel\_yaqubi@yahoo.es}

\thanks{$^*$Corresponding author: Daniel Yaqubi (\texttt{yaqubi@tjamcaas.ac.ir}).}
\author{Madjid Mirzavaziri}
\address{Department of Pure Mathematics, Ferdowsi University of Mashhad, Mashhad, Iran}
\email{mirzavaziri@um.ac.ir, mirzavaziri@gmail.com}

\date{}

\begin{document}
\begin{abstract}
Let $A_k(n)$ and $B_k(n)$ be the numbers of ordered $k$-tuples of partitions with designated summands, respectively with designated summands and all parts odd, of total size $n$. Using an exact description of the generating function of designated partitions after a $3$-dissection, we classify for every $e$ the values of $k$ for which $A_k(3n+2)$ is divisible by $3^e$: this happens exactly when $k\equiv0$ or $-5\pmod{3^e}$ if $e\le4$, and when $3^e\mid k$ if $e\ge5$. We obtain the corresponding exact $3$-adic and several $2$-adic results for $B_k$ and for $A_k$ along the progressions $4n+3$, a recursion for the classes $3^b(3n+s)$, and for $k=2$ the congruences $A_2(12n+10)\equiv 27f_1^4f_6^3/f_2\pmod{81}$ in generating-function form. We relate the numbers $A_k(n)$ modulo $3$ and $9$ to overpartitions, sums of squares and triangular numbers, give group-action proofs of some of the congruences, and close with several conjectures supported by computation.
\end{abstract}
\maketitle

\section{Introduction}\label{sec:intro}
Since Ramanujan’s famous partition congruences $p(5n+4)\equiv0\pmod 5$, $p(7n+5)\equiv0\pmod 7$, and $p(11n+6)\equiv0\pmod{11}$, the arithmetic properties of restricted partition families have been a central focus in number theory. While traditional methods often rely on generating functions, modular forms, or combinatorial statistics like Dyson’s rank and the Andrews–Garvan crank, looking at $k$-tuples (multipartitions) reveals deeper structures whose properties depend heavily on the $p$-adic valuation of $k$.

\begin{figure}[htbp]
\centering
\small
\begin{tikzpicture}[
    x=0.35cm, y=0.35cm,
    cell/.style={draw=blue!50!black, fill=blue!8!white, line width=0.6pt, rounded corners=0.5pt},
    dot/.style={fill=red!75!black}
]
    \draw[cell] (0.0,0.0) rectangle (1.0,-1.0);
    \draw[cell] (1.0,0.0) rectangle (2.0,-1.0);
    \draw[cell] (2.0,0.0) rectangle (3.0,-1.0);
    \draw[cell] (3.0,0.0) rectangle (4.0,-1.0);
    \fill[dot] (0.5,-0.5) circle (0.22);

    \draw[cell] (5.2,0.0) rectangle (6.2,-1.0);
    \draw[cell] (6.2,0.0) rectangle (7.2,-1.0);
    \draw[cell] (7.2,0.0) rectangle (8.2,-1.0);
    \fill[dot] (5.7,-0.5) circle (0.22);
    \draw[cell] (5.2,-1.0) rectangle (6.2,-2.0);
    \fill[dot] (5.7,-1.5) circle (0.22);

    \draw[cell] (10.4,0.0) rectangle (11.4,-1.0);
    \draw[cell] (11.4,0.0) rectangle (12.4,-1.0);
    \fill[dot] (10.9,-0.5) circle (0.22);
    \draw[cell] (10.4,-1.0) rectangle (11.4,-2.0);
    \draw[cell] (11.4,-1.0) rectangle (12.4,-2.0);

    \draw[cell] (15.6,0.0) rectangle (16.6,-1.0);
    \draw[cell] (16.6,0.0) rectangle (17.6,-1.0);
    \draw[cell] (15.6,-1.0) rectangle (16.6,-2.0);
    \draw[cell] (16.6,-1.0) rectangle (17.6,-2.0);
    \fill[dot] (16.1,-1.5) circle (0.22);

    \draw[cell] (20.8,0.0) rectangle (21.8,-1.0);
    \draw[cell] (21.8,0.0) rectangle (22.8,-1.0);
    \fill[dot] (21.3,-0.5) circle (0.22);
    \draw[cell] (20.8,-1.0) rectangle (21.8,-2.0);
    \fill[dot] (21.3,-1.5) circle (0.22);
    \draw[cell] (20.8,-2.0) rectangle (21.8,-3.0);

    \draw[cell] (0.0,-6.2) rectangle (1.0,-7.2);
    \draw[cell] (1.0,-6.2) rectangle (2.0,-7.2);
    \fill[dot] (0.5,-6.7) circle (0.22);
    \draw[cell] (0.0,-7.2) rectangle (1.0,-8.2);
    \draw[cell] (0.0,-8.2) rectangle (1.0,-9.2);
    \fill[dot] (0.5,-8.7) circle (0.22);

    \draw[cell] (5.2,-6.2) rectangle (6.2,-7.2);
    \fill[dot] (5.7,-6.7) circle (0.22);
    \draw[cell] (5.2,-7.2) rectangle (6.2,-8.2);
    \draw[cell] (5.2,-8.2) rectangle (6.2,-9.2);
    \draw[cell] (5.2,-9.2) rectangle (6.2,-10.2);

    \draw[cell] (10.4,-6.2) rectangle (11.4,-7.2);
    \draw[cell] (10.4,-7.2) rectangle (11.4,-8.2);
    \fill[dot] (10.9,-7.7) circle (0.22);
    \draw[cell] (10.4,-8.2) rectangle (11.4,-9.2);
    \draw[cell] (10.4,-9.2) rectangle (11.4,-10.2);

    \draw[cell] (15.6,-6.2) rectangle (16.6,-7.2);
    \draw[cell] (15.6,-7.2) rectangle (16.6,-8.2);
    \draw[cell] (15.6,-8.2) rectangle (16.6,-9.2);
    \fill[dot] (16.1,-8.7) circle (0.22);
    \draw[cell] (15.6,-9.2) rectangle (16.6,-10.2);

    \draw[cell] (20.8,-6.2) rectangle (21.8,-7.2);
    \draw[cell] (20.8,-7.2) rectangle (21.8,-8.2);
    \draw[cell] (20.8,-8.2) rectangle (21.8,-9.2);
    \draw[cell] (20.8,-9.2) rectangle (21.8,-10.2);
    \fill[dot] (21.3,-9.7) circle (0.22);
\end{tikzpicture}
\caption{The ten partitions with designated summands of $n = 4$ ($\mathrm{PD}(4) = 10$) represented via Young diagrams with designated dotted cells.}
\label{fig:young}
\end{figure}
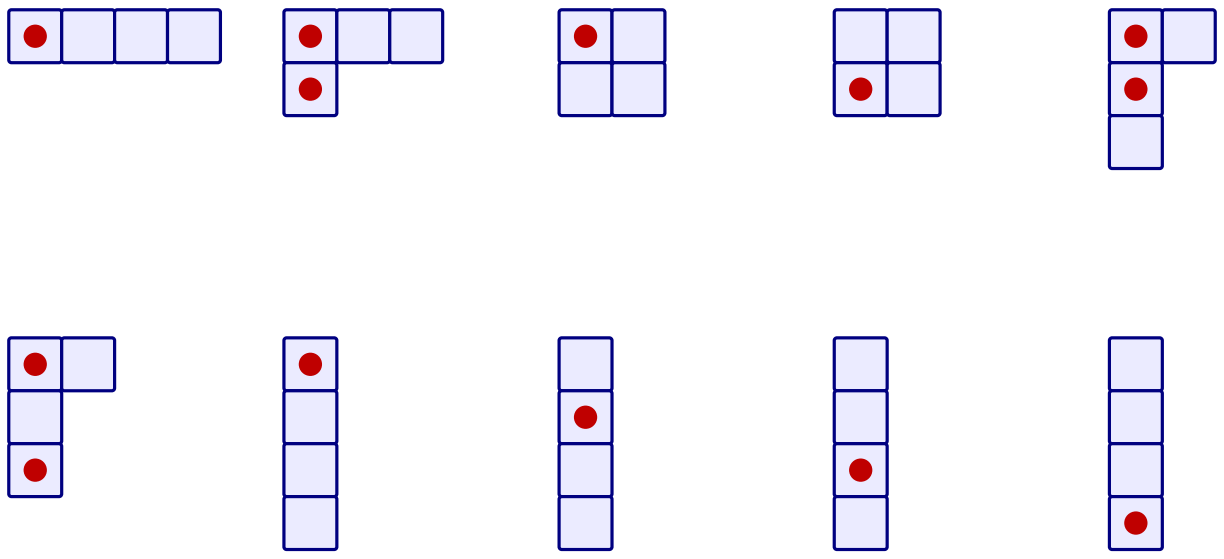

A prominent family in this area is partitions with designated summands, $\mathrm{PD}(n)$, introduced by Andrews, Lewis, and Lovejoy, where exactly one part of each size is marked. Their generating function is given by $\sum \mathrm{PD}(n)q^n = f_6/(f_1f_2f_3)$, where $f_k := \prod_{n=1}^\infty (1-q^{kn})$ denotes the standard eta-quotient. Figure~\ref{fig:young} illustrates this combinatorial structure through Young diagrams for $n=4$, where every block of equal rows contains a designated marked cell. Over time, these objects have been extensively studied using $3$-dissections, cranks, and various congruence families. Similarly, restricting parts to be odd gives $\mathrm{PDO}(n)$, explored by Baruah and Ojah and later refined by Fu and Sellers.

For bipartitions with designated summands, $A_2(n) = \mathrm{PD}_{-2}(n)$, Naika and Shivashankar along with Hao and Shen proved several deep congruences modulo powers of $2$ and $3$. Motivated by these works, we study the families $A_k(n) = \mathrm{PD}_{-k}(n)$ and $B_k(n) = \mathrm{PDO}_{-k}(n)$, which count ordered $k$-tuples of designated partitions and odd-part designated partitions. Our approach starts with the $3$-dissection of $f_1f_2$ using Hirschhorn's classical methods, writing the main generating function as $F = f_6/(f_1f_2f_3) = \Theta P(u)$ with a quartic polynomial $P$ and explicit eta-quotients $\Theta$ and $u$. This turns questions about $F^k \pmod{3^e}$ into finite checks on the coefficients of $P(u)^k$.

Throughout this paper, we give a complete classification for every $e \ge 1$ of the parameters $k$ where $A_k(3n+2) \equiv 0 \pmod{3^e}$ for all $n$, and we provide level-wise recursions. We also prove exact formulas and new congruences for $k=2$ that improve earlier results, such as $A_2(24n+23) \equiv 0 \pmod{64}$ and $A_2(18n+15) \equiv 0 \pmod{32}$, alongside general congruences for $A_k(4n+3)$ and the exact $3$-adic valuation $v_3(k)+1$ for $B_k(3n+2)$ when $3 \mid k$. The rest of the paper explores combinatorial interpretations using Young diagrams, links to overpartitions and sums of squares, group-action proofs, and computational conjectures, all using elementary tools like Frobenius congruences and Lucas' theorem.
\section{Preliminaries and the structure of $F$}\label{sec:prelim}
In this section, we review basic notation. Let $f_m = (q^m;q^m)_\infty = \prod_{j\ge1}(1-q^{mj})$ for $m \ge 1$. Products of the form $\prod_\delta f_\delta^{r_\delta}$ are called eta-quotients, which we use mainly for notation and dissection identities. We also frequently use standard product relations such as $(-q;q)_\infty = f_2/f_1$, $(q;q^2)_\infty = f_1/f_2$, and $(-q;q^2)_\infty = f_2^2/(f_1f_4)$.

Ramanujan's theta functions are defined by $\varphi(q) = \sum_{m\in\mathbb Z}q^{m^2} = f_2^5/(f_1^2f_4^2)$ and $\psi(q) = \sum_{m\ge0}q^{m(m+1)/2} = f_2^2/f_1$. We also use the Borwein cubic theta functions $a(q) = \sum_{m,n}q^{m^2+mn+n^2}$ and $b(q) = f_1^3/f_3 = \frac12(3a(q^3)-a(q))$ \cite{BBG}. For comparison, we consider overpartitions $\overline{p}(n)$ with generating function $f_2/f_1^2 = 1/\varphi(-q)$, and cubic partitions $a(n)$ with generating function $1/(f_1f_2)$. For any partition family $X$, we write $X_{-k}(n)$ to denote the number of $k$-tuples of total size $n$, whose generating function is the $k$-th power of the original one.

Recall that a partition with designated summands is defined by marking exactly one part of each size, which corresponds on its Young diagram to having exactly one row marked with a dot in each block of equal rows (Figure~\ref{fig:young}). Since a block of $m$ equal rows can be marked in $m$ ways, the generating function is given by
\[
\sum_{n\ge0}\mathrm{PD}(n)q^n = \prod_{t\ge1}\left(1+\sum_{m\ge1}mq^{tm}\right) = \prod_{t\ge1}\left(1+\frac{q^t}{(1-q^t)^2}\right) = \frac{f_6}{f_1f_2f_3}.
\]
Similarly, when all parts are restricted to be odd, we obtain $\mathrm{PDO}(n)$, whose generating function is $\sum \mathrm{PDO}(n)q^n = f_4f_6^2/(f_1f_3f_{12})$.
\begin{definition}
For $k\ge1$, let $A_k(n)=\mathrm{PD}_{-k}(n)$ and $B_k(n)=\mathrm{PDO}_{-k}(n)$, defined by their generating functions:
\begin{align*}
\sum_{n\ge0}A_k(n)q^n&=F^k,\qquad F=\frac{f_6}{f_1f_2f_3}=1+q+3q^2+5q^3+10q^4+\cdots,\\
\sum_{n\ge0}B_k(n)q^n&=G^k,\qquad G=\frac{f_4f_6^2}{f_1f_3f_{12}}=1+q+2q^2+4q^3+5q^4+\cdots.
\end{align*}
\end{definition}

The \emph{$m$-dissection} of a series $\sum c_nq^n$ splits it into $\sum_{r=0}^{m-1}q^rD_r(q^m)$, where $D_r(q^m)=\sum_{n\equiv r\,(m)}c_nq^{n-r}$. For $m=2$, this gives the even and odd parts. We use the classical $2$- and $3$-dissections of eta-quotients due to Hirschhorn \cite{Hir}, which we cite as needed. 

In addition to dissections, we rely on two standard tools: Lucas' theorem, which states that for prime $p$ with base-$p$ expansions $n=\sum n_ip^i$ and $m=\sum m_ip^i$, we have $\binom{n}{m}\equiv\prod_i\binom{n_i}{m_i}\pmod p$; and the standard Frobenius congruence.

\begin{lemma}[\cite{Hir}]\label{lem:frob}
For any prime $p$, integer $j\ge1$, and series $H$ with integer coefficients,
\[
H(q)^{p^j}\equiv H(q^p)^{p^{j-1}}\pmod{p^j}.\]
\end{lemma}

A first direct consequence is that the divisibility properties of terms in the arithmetic progressions $3n+r$ are periodic with respect to $k$.

\begin{lemma}\label{lem:period}
Let $r\in\{0,1,2\}$ and $e\ge1$. Then $A_k(3n+r)\equiv0\pmod{3^e}$ for all $n$ if and only if $A_{k+3^e}(3n+r)\equiv0\pmod{3^e}$ for all $n$. The same property holds for $B_k$.
\end{lemma}

\begin{proof}
By Lemma~\ref{lem:frob}, we have $F^{3^e}\equiv U(q^3)\pmod{3^e}$, where $U(q^3)=F(q^3)^{3^{e-1}}$ with $U(0)=1$. This gives $F^{k+3^e}\equiv U(q^3)F^k$. Writing $U(q^3)=\sum_iu_iq^{3i}$, it follows that
\[
A_{k+3^e}(3n+r)\equiv\sum_{i\le n}u_iA_k(3(n-i)+r)\pmod{3^e}.
\]
If all $A_k(3m+r)$ are divisible by $3^e$, the left side is clearly divisible as well. Conversely, since $u_0=1$, induction on $n$ shows that $A_k(3n+r)$ is also divisible by $3^e$.
\end{proof}
Our main tool for all subsequent results is an explicit description of the series $F$.

\begin{proposition}\label{prop:struct}
Let
\[
c=\frac{f_6f_9^4}{f_3f_{18}^2},\qquad \rho=\frac{f_3f_{18}^3}{f_6f_9^3},\qquad u=q\rho,\qquad \Lambda=\frac{f_{18}^2}{f_9^4}.
\]
Then $c,\rho,\Lambda\in\mathbb Z[[q^3]]$ have constant term $1$, and
\[
f_1f_2=c\,(1-u-2u^2),\qquad F=\frac{\Lambda}{1-u-2u^2}=\Theta\,P(u),
\]
where $P(u)=1+u+3u^2-2u^3+4u^4=(1-u+u^2)(1+2u+4u^2)$ and $\Theta=\Lambda/\bigl((1+u^3)(1-8u^3)\bigr)\in\mathbb Z[[q^3]]$ has constant term $1$.
\end{proposition}
\begin{proof}
By Hirschhorn's $3$-dissection (see \cite[Lemma~2.6]{NHPV}), we have \[f_1f_2=\frac{f_6f_9^4}{f_3f_{18}^2}-q\,f_9f_{18}-2q^2\frac{f_3f_{18}^4}{f_6f_9^2}.\]Direct cancellation gives $f_9f_{18}=c\rho$ and $f_3f_{18}^4/(f_6f_9^2)=c\rho^2$, which leads to $f_1f_2=c(1-u-2u^2)$. Since $f_6/(f_3c)=f_{18}^2/f_9^4=\Lambda$, we get $F=f_6/(f_3f_1f_2)=\Lambda/(1-u-2u^2)$. Finally, since $1-u-2u^2=(1+u)(1-2u)$, multiplying both the numerator and denominator by $(1-u+u^2)(1+2u+4u^2)$ turns the denominator into $(1+u^3)(1-8u^3)$.
\end{proof}

While the formula $F=\Lambda/(1-u-2u^2)$ is equivalent to the $3$-dissection from \cite{CJJS}, our new polynomial form $F=\Theta P(u)$ is the key insight: it reduces any question about $F^k$ modulo $3^e$ to a finite coefficient check (Lemma~\ref{lem:class}).

\begin{lemma}\label{lem:class}
Let $c_j(k)=[u^j]P(u)^k$. For any $r\in\{0,1,2\}$ and $M\ge1$,
\[
A_k(3n+r)\equiv0\pmod M\ \text{ for all }n\iff c_j(k)\equiv0\pmod M\ \text{ for all }j\equiv r\pmod 3.
\]
\end{lemma}
\begin{proof}
We can write $F^k = \Theta^k \sum_j c_j(k) q^j \rho^j$. Since $\Theta$ and $\rho$ belong to $\mathbb Z[[q^3]]$, the terms of $F^k$ with exponents congruent to $r \pmod 3$ are given by $q^r \Theta^k \rho^r \sum_{i\ge0} c_{r+3i}(k)(q^3\rho^3)^i$. Here, $\Theta^k \rho^r$ is a unit in $\mathbb Z[[q^3]]$, and $(q^3\rho^3)^i$ has exact order $3i$ with leading coefficient $1$. Thus, all coefficients of $q^{3n+r}$ vanish modulo $M$ if and only if every $c_{r+3i}(k) \equiv 0 \pmod M$ (which follows by inspecting the smallest $i$ where $c_{r+3i}$ is non-zero).
\end{proof}

We will also need the explicit eta-quotient form of $\Theta$ and a factorization of $(1+u^3)(1-8u^3)$.
\begin{lemma}\label{lem:EW}
Let $w = u^3|_{q^3 \to x} = x \, f_1^3 f_6^9 / (f_2^3 f_3^9)$, where $f_m = f_m(x)$. Then,
\[
1+w = \frac{f_2^5 f_6}{f_1 f_3^5}, \qquad 1-8w = \frac{f_1^8 f_6^4}{f_2^4 f_3^8}, \qquad \text{and consequently} \qquad \Theta = \frac{f_9^9}{f_3^7 f_6 f_{18}^3}.
\]
\end{lemma}
\begin{proof}
Both sides of each identity are modular functions of weight $0$ on $\Gamma_0(6)$ and are analytic (holomorphic) on the upper half-plane. By the valence formula (since the index of $\Gamma_0(6)$ is $12$), the difference between the two sides is uniquely determined by a finite number of coefficients in its $q$-expansion, which we can check to be equal. 

To prove the final expression for $\Theta$, we multiply the first two identities to get $(1+w)(1-8w) = f_1^7 f_2 f_6^5 / f_3^{13}$. In terms of $q$, this becomes $(1+u^3)(1-8u^3) = f_3^7 f_6 f_{18}^5 / f_9^{13}$. (Note: This product can also be proved elementarily without using modularity by starting from $f_1f_2=c(1-u-2u^2)$ and multiplying over the three cube roots of unity $q \mapsto \omega^j q$). Finally, dividing $\Lambda = f_{18}^2 / f_9^4$ by this product yields the formula for $\Theta$.
\end{proof}

\section{Congruences for $A_k$}\label{sec:A}

This section contains our main results. We begin by examining the $3$-adic behavior of $A_k$: we completely classify the arithmetic progression $3n+2$ and describe the classes $3^b(3n+s)$ modulo $3$. We then turn to powers of $2$ and the bipartition case $k=2$, where an exact formula for the class $3n+1$ lifts a known congruence modulo $9$ to moduli $27$ and $81$. We first state the elementary divisibility that follows directly from Lemma~\ref{lem:frob}.

\begin{theorem}\label{thm:frob}
If $p$ is a prime, $p^j\mid k$, and $1\le r\le p-1$, then $A_k(pn+r)\equiv B_k(pn+r)\equiv0\pmod{p^j}$.
\end{theorem}
\begin{proof}
By Lemma~\ref{lem:frob}, $F^{k}=(F^{p^j})^{k/p^j}\equiv\bigl(F(q^p)^{p^{j-1}}\bigr)^{k/p^j}\pmod{p^j}$, which is a series in $q^p$.
\end{proof}

This represents the standard argument for $k$-tuples. The more interesting progressions are those it does not reach, which form the focus of the rest of the paper.

The progression $3n+2$ is more delicate and is governed by the negative power $F^{-5}$. We write $A_{-5}(n)=[q^n]F^{-5}=[q^n]\,(f_1f_2f_3)^5/f_6^5$, even though this coefficient does not count combinatorial objects.

\begin{proposition}\label{prop:neg5}
$\displaystyle\sum_{n\ge0}A_{-5}(3n+2)\,x^n=-81\,x\,\frac{f_1^5f_3^5f_6^5}{f_2^5}$.
\end{proposition}
\begin{proof}
By Proposition~\ref{prop:struct}, $F^{-5}=\Lambda^{-5}(1-u-2u^2)^5$, and
\[
(1-u-2u^2)^5=1-5u+30u^3-15u^4-81u^5+30u^6+120u^7-80u^9-32u^{10}.
\]
The only term with an exponent congruent to $2\pmod 3$ is $-81u^5$ (since the coefficients of $u^2$ and $u^8$ vanish). Because $\Lambda$ and $\rho$ belong to $\mathbb Z[[q^3]]$, the part of $F^{-5}$ with exponents congruent to $2\pmod3$ is $-81\,q^5\rho^5\Lambda^{-5}$. Dividing by $q^2$, setting $x=q^3$, and using the substitutions $\rho\mapsto f_1f_6^3/(f_2f_3^3)$ and $\Lambda\mapsto f_6^2/f_3^4$ yields the result.
\end{proof}

\begin{theorem}\label{thm:class3}
Let $k,e\ge1$. Then $A_k(3n+2)\equiv0\pmod{3^e}$ for all $n\ge0$ if and only if
\[
\begin{cases}
k\equiv0\ \text{or}\ k\equiv-5\pmod{3^e}, & \text{if }e\le4, \\ 
k\equiv0\pmod{3^e}, & \text{if }e\ge5.
\end{cases}
\]
\end{theorem}
\begin{proof}
Since $A_k(2)=3k+\binom k2=k(k+5)/2$ and $\gcd(k,k+5)\mid5$, a necessary condition is that $3^e$ divides either $k$ or $k+5$. If $3^e\mid k$, the congruence holds by Theorem~\ref{thm:frob}. Let $k=3^em-5$ with $m\ge1$. Then $F^k=(F^{3^e})^mF^{-5}$, and by Lemma~\ref{lem:frob}, $(F^{3^e})^m\equiv V(q^3)\pmod{3^e}$ for some $V\in\mathbb Z[[q^3]]$ with $V(0)=1$. Thus, by Proposition~\ref{prop:neg5},
\[
\sum_nA_k(3n+2)x^n\equiv-81\,x\,V(x)\,\frac{f_1^5f_3^5f_6^5}{f_2^5}\pmod{3^e}.
\]
For $e\le4$, the right-hand side vanishes modulo $3^e$. For $e\ge5$, the coefficient of $x$ is $-81\not\equiv0\pmod{3^5}$, meaning $3^e\nmid A_k(5)$.
\end{proof}

\begin{remark}
For $e=1,2$, the theorem states that
\begin{align*}
A_k(3n+2)\equiv0\pmod 3&\iff k\equiv0,1\pmod 3,\\
A_k(3n+2)\equiv0\pmod 9&\iff k\equiv0,4\pmod 9.
\end{align*}
For $e=3,4$, the exceptional residues are $k\equiv22\pmod{27}$ and $k\equiv76\pmod{81}$. In particular, modulo $3^5$, the progression $3n+2$ yields no additional congruences beyond $3^5\mid k$. For every $k\equiv-5\pmod{81}$, the highest power of $3$ dividing all $A_k(3n+2)$ is exactly $3^4$: if $243\mid k+5$, then $A_k(5)\equiv-81\pmod{243}$; otherwise, $v_3(A_k(2))=4$. The case $k=1, e=1$ corresponds to the classical congruence of Andrews, Lewis, and Lovejoy \cite{ALL}.
\end{remark}

For $e\le2$, there is an alternative proof that avoids Hirschhorn's dissection of $f_1f_2$ and explicitly connects the numbers $A_k(n)$ to overpartitions, sums of squares, and the divisor function of the hexagonal lattice. Let $\overline p_{-m}(n)$ denote the number of $m$-tuples of overpartitions of total size $n$, such that the generating function is given by $\sum\overline p_{-m}(n)q^n=\varphi(-q)^{-m}$.
\begin{proposition}\label{prop:dict}
We have $F \equiv \varphi(-q)^{-2}$ and $G \equiv \varphi(q)/\varphi(-q) \pmod 3$. Consequently, for all $k \ge 1$ and $n \ge 0$,
\[
A_k(n) \equiv \overline{p}_{-2k}(n) \pmod 3, \qquad
B_k(n) \equiv \sum_{j} r_k(j) \, \overline{p}_{-k}(n-j) \pmod 3,
\]
where $r_k(j)$ denotes the number of representations of $j$ as an ordered sum of $k$ squares.
\end{proposition}

\begin{proof}
Working modulo $3$, we use the standard basic relations $f_3 \equiv f_1^3$, $f_6 \equiv f_2^3$, and $f_{12} \equiv f_4^3$. 

First, substituting these relations into the definition of $F$ gives:
\[
F = \frac{f_6}{f_1 f_2 f_3} \equiv \frac{f_2^3}{f_1 f_2 f_1^3} = \frac{f_2^2}{f_1^4} = \left(\frac{f_2}{f_1^2}\right)^2 \pmod 3.
\]
Since the generating function for overpartitions satisfies $\sum \overline{p}(n)q^n = f_2/f_1^2 = \varphi(-q)^{-1}$, it immediately follows that $F \equiv \varphi(-q)^{-2} \pmod 3$.

Second, for $G$, we start with the classical theta-function identity $\frac{\varphi(q)}{\varphi(-q)} = \frac{f_2^6}{f_1^4 f_4^2}$. To prove that $G \equiv \frac{f_2^6}{f_1^4 f_4^2} \pmod 3$, we multiply $G$ by the reciprocal expression and simplify modulo $3$:
\[
G \cdot \frac{f_1^4 f_4^2}{f_2^6} = \left(\frac{f_4 f_6^2}{f_1 f_3 f_{12}}\right) \cdot \frac{f_1^4 f_4^2}{f_2^6} = \frac{f_1^3 f_4^3 f_6^2}{f_1^3 f_4^3 f_2^6}.
\]
Using $f_6 \equiv f_2^3 \pmod 3$, the numerator becomes $f_1^3 f_4^3 (f_2^3)^2 = f_1^3 f_4^3 f_2^6$, which exactly matches the denominator. Thus, the product is congruent to $1 \pmod 3$, proving that $G \equiv \frac{f_2^6}{f_1^4 f_4^2} = \frac{\varphi(q)}{\varphi(-q)} \pmod 3$.

Finally, raising these congruences to the power $k$ yields the desired coefficient relations for $A_k(n) = [q^n]F^k$ and $B_k(n) = [q^n]G^k$ in terms of overpartition tuples and sums of squares.
\end{proof}
\begin{corollary}\label{cor:dict}
For $j \ge 0$, $F^{3j+1} \equiv \varphi(-q) \, \varphi(-q^3)^{-2j-1}$ and $F^{3j+2} \equiv \varphi(-q)^2 \varphi(-q^3)^{-2j-2} \pmod 3$. This clearly shows the dependence on $k \bmod 3$: $\varphi(-q)$ has no exponents congruent to $2 \pmod 3$, while $\varphi(-q)^2$ has coefficient $4$ at $q^2$.
\end{corollary}
\begin{proof}
This follows directly from $\varphi(-q)^{-3} \equiv \varphi(-q^3)^{-1}$ (Lemma~\ref{lem:frob}).
\end{proof}
\begin{proposition}\label{prop:dict9}
Let $T(q) = \sum_{n\ge1} S(n)q^n - \sum_{n\ge1} S(n)q^{2n}$. Then
\[
F = \varphi(-q)^{-2} \, \frac{b(q)}{b(q^2)}, \qquad
F^k \equiv \varphi(-q)^{-2k} \bigl(1 - 3k \, T(q)\bigr) \pmod 9.
\]
Equivalently, $A_k(n) \equiv \overline{p}_{-2k}(n) - 3k \sum_j t(j) \, \overline{p}_{-2k}(n-j) \pmod 9$, where $t(j) = S(j) - S(j/2)$ and $S(j/2) = 0$ for odd $j$.
\end{proposition}

\begin{proof}
We start with the exact identity $F \varphi(-q)^2 = f_1^3 f_6 / (f_2^3 f_3) = b(q)/b(q^2)$. 
Since $b(q) \equiv 1 - 3\sum S(n)q^n \pmod 9$, working modulo $9$ gives:
\[
\frac{b(q)}{b(q^2)} \equiv \frac{1 - 3\sum S(n)q^n}{1 - 3\sum S(n)q^{2n}} \equiv (1 - 3\sum S(n)q^n)(1 + 3\sum S(n)q^{2n}) = 1 - 3T(q) \pmod 9.
\]
Raising both sides to the power $k$ yields the desired congruence for $F^k$.
\end{proof}
In other words, modulo $9$, the numbers $A_k(n)$ are determined by overpartition tuples combined with the divisor function of the hexagonal lattice.

\begin{proof}[Second proof of Theorem~\ref{thm:class3} for $e \le 2$]
For $e=1$, the statement follows from Corollary~\ref{cor:dict} together with the coefficient $4$ at $q^2$ of $\varphi(-q)^2$ when $k \equiv 2 \pmod 3$. For $e=2$, necessity follows from $A_k(2) = k(k+5)/2$, while the case $9 \mid k$ is covered by Theorem~\ref{thm:frob}. By Lemma~\ref{lem:period}, it remains to check $k=4$. Writing $\varphi = \varphi(-q)$ and $\psi = \varphi(-q^3)$, this proof relies entirely on classical theta-function dissections. 

By Proposition~\ref{prop:dict9}, $F^4 \equiv \varphi^{-8}(1-12T) \equiv \varphi^{-8}(1-3T) \pmod 9$. Since $\varphi^9 \equiv \psi^3 \pmod 9$ (Lemma~\ref{lem:frob}), we have $\varphi^{-8} = \varphi \cdot \varphi^{-9} \equiv \varphi \, \psi^{-3}$, which leads to
\[
F^4 \equiv \psi^{-3} \, \varphi \, (1-3T) \pmod 9,
\]
where $\psi^{-3}$ is a unit series in $q^3$. Therefore, it suffices to show that all coefficients of $q^{3n+2}$ in $\varphi(1-3T)$ are divisible by $9$. Since $\varphi = \sum (-1)^m q^{m^2}$ has no exponents congruent to $2 \pmod 3$, this reduces to showing that
\[
[q^{3n+2}] \, \varphi \, T \equiv 0 \pmod 3.
\]
We apply a $3$-dissection to the components. Let $S(q) = \sum_{n\ge1} S(n)q^n$. Using Borwein's dissections $a(q) = a(q^3) + 6q \, f_9^3/f_3$ and $a = 1 + 6S$, we obtain
\[
S(q) = S(q^3) + q \frac{f_9^3}{f_3}, \quad \text{implying that } S(3n+2) = 0.
\]
Replacing $q$ by $q^2$ and splitting according to $n \bmod 3$ gives
\[
S(q^2) = S(q^6) + q^2 \frac{f_{18}^3}{f_6}.
\]
Consequently, $T = \bigl[S(q^3) - S(q^6)\bigr] + q \, f_9^3/f_3 - q^2 f_{18}^3/f_6$, with the three summands falling into classes $0, 1, and 2$ modulo $3$, respectively. From Hirschhorn's dissection $\varphi(-q) = \varphi(-q^9) - 2q \, f_3 f_{18}^2 / (f_6 f_9)$, the class-$0$ part is $f_9^2 / f_{18}$ and the class-$1$ part is $-2q \, f_3 f_{18}^2 / (f_6 f_9)$. The class-$2$ part of $\varphi T$ is therefore given by
\[
-\frac{f_9^2}{f_{18}} \cdot q^2 \frac{f_{18}^3}{f_6} \;-\; 2q \frac{f_3 f_{18}^2}{f_6 f_9} \cdot q \frac{f_9^3}{f_3} 
= -3q^2 \frac{f_9^2 f_{18}^2}{f_6} \equiv 0 \pmod 3.
\]
This completes the proof for $k \equiv 4 \pmod 9$.
\end{proof}
We next record the classes $9n+3$ and $9n+6$ and, more generally, $3^b(3n+s)$ modulo $3$.

\begin{theorem}\label{thm:mod9cls}
For every $k\ge0$:
\begin{enumerate}
\item[(a)] $A_k(9n+3)\equiv0\pmod 3$ for all $n$ if and only if $k\equiv0,4,8\pmod 9$;
\item[(b)] $A_k(9n+6)\equiv0\pmod 3$ for all $n$ if and only if $k\equiv0,2,3,4,7,8\pmod 9$.
\end{enumerate}
\end{theorem}
\begin{proof}
All congruences take place in $\mathbb F_3[[x]]$. Let $c,\rho,u=x\rho,\Lambda$ be as in Proposition~\ref{prop:struct}, with $q$ renamed $x$.

\emph{Step 1: $F$ modulo $3$.} Since $-2\equiv1\pmod 3$ we have $1-u-2u^2\equiv(1+u)^2$, and therefore
\[
f_1f_2\equiv c\,(1+u)^2,\qquad
F=\frac{\Lambda}{1-u-2u^2}\equiv\frac{\Lambda}{(1+u)^2}=\frac{\Lambda\,(1+u)}{(1+u)^3}\equiv\kappa\,(1+u),
\]
where $\kappa=f_{18}^2/\bigl(f_9^4(1+u^3)\bigr)\in\mathbb F_3[[x^3]]$.

\emph{Step 2: exponent classes.} Put $R=\mathbb F_3[[x^3]]$. Since $\rho\in R$ is a unit and the orders of the series $u^{3i}$ are strictly increasing, a series $\sum_ic_iu^{3i+s}$ with $c_i\in\mathbb F_3$ vanishes if and only if all $c_i$ do. In this sense the classes of exponents of $x$ modulo $3$ correspond to the classes of exponents of $u$ modulo $3$.

\emph{Step 3: the part $3m$.} Write $k=3t+r$ with $0\le r\le2$. Then $F^k\equiv\kappa^k\sum_j\binom kj u^j$, and by Lucas' theorem $\binom k{3i}\equiv\binom ti$. So the part of $F^k$ with exponents divisible by $3$ is $\kappa^k\sum_i\binom tiu^{3i}$. Writing $y=x^3$ and using $u(x)^3\equiv u(y)$ we obtain
\[
\sum_{m\ge0}A_k(3m)\,y^m\equiv\Bigl(\frac{f_6^2}{f_3^4}\Bigr)^{k}(1+u)^{t-k},\qquad u=u(y)=y\rho(y),
\]
where now $f_j=f_j(y)$ and $f_6^2/f_3^4$ is a unit of $\mathbb F_3[[y^3]]$.

\emph{Step 4: the vanishing criterion.} Choose $e$ with $N=3e+t-k\ge0$. Then $(1+u)^{t-k}=(1+u)^N/(1+u^3)^e$, and $(1+u^3)^e$ is a unit of $\mathbb F_3[[y^3]]$. By Step~2, the coefficients of $y^{3n+s}$ ($s=1,2$) vanish for all $n$ if and only if $\binom N{3i+s}\equiv0$ for all $i$. By Lucas' theorem
\[
\binom N{3i+s}\equiv\binom{N_0}{s}\binom{\lfloor N/3\rfloor}{i},\qquad N_0=N\bmod 3,
\]
so the condition is $N_0<s$, i.e.\ $(t-k)\bmod 3<s$. As $t-k=-2t-r\equiv t-r\pmod 3$, this reads $(t-r)\bmod 3<s$.

\emph{Step 5: the residues.} Running through the nine residues $k\bmod 9=3(t\bmod 3)+r$, the condition $(t-r)\bmod 3=0$ (for $s=1$) holds exactly for $k\equiv0,4,8\pmod 9$, and the condition $(t-r)\bmod 3\in\{0,1\}$ (for $s=2$) holds exactly for $k\equiv0,2,3,4,7,8\pmod 9$. Finally $A_k(9n+3)=A_k(3(3n+1))$ and $A_k(9n+6)=A_k(3(3n+2))$.
\end{proof}

\begin{theorem}\label{thm:levels}
Let integers $\lambda_a$ be defined by $\lambda_0=0$, $\lambda_1=-4$, and $\lambda_{a+1}=-\lambda_a+4\lambda_{a-1}$ (giving $\lambda_a = -4, 4, -20, 36, -116, \dots$), and for $k\ge0$, let $M_1=\lfloor k/3\rfloor-k$ and $M_{a+1}=\lfloor M_a/3\rfloor+2\lambda_a k$. Then, for all $b\ge1$, $s\in\{1,2\}$, and $k\ge0$,
\[
A_k\bigl(3^b(3n+s)\bigr)\equiv0\pmod 3\ \text{ for all }n\ge0\iff M_b(k)\bmod 3<s.
\]
\end{theorem}

\begin{proof}
The argument iterates the method of Theorem~\ref{thm:mod9cls} using the notation $c, \rho, u, R$. Let $\Phi_a(z)=\sum_{m\ge0}A_k(3^am)z^m\in\mathbb F_3[[z]]$. We first establish by induction that for every $a\ge1$,
\[
\Phi_a(z)\equiv W_a(z)\,\bigl(1+u(z)\bigr)^{M_a},
\]
where $W_a = f_3^{p_3}f_6^{p_6}f_9^{p_9}f_{18}^{p_{18}}\in\mathbb F_3[[z^3]]$ with exponent vector $p^{(a)}=(p_3,p_6,p_9,p_{18})$. For $a=1$, this holds by Theorem~\ref{thm:mod9cls} with $p^{(1)}=k\,(-4,2,0,0)$. Assuming the relation holds for some $a$ in terms of the variable $y$, we choose $e$ such that $N = M_a + 3e \ge 0$, allowing us to write $(1+u(y))^{M_a} = (1+u(y))^N(1+u(y)^3)^{-e}$. By Lucas' theorem, the part of $(1+u)^N$ with exponents congruent to $0 \pmod 3$ is $(1+u^3)^{\lfloor N/3\rfloor}$. Since $W_a(y)\in\mathbb F_3[[y^3]]$ and $u(y)^3=u(y^3)$, setting $z=y^3$ yields
\[
\Phi_{a+1}(z)\equiv W_a(y)\,\bigl(1+u(z)\bigr)^{\lfloor M_a/3\rfloor}.
\]
Using $f_{3j}(y)=f_j(z)$, $W_a(y)$ evaluates to $f_1^{p_3}f_2^{p_6}f_3^{p_9}f_6^{p_{18}}$ at $z$. By induction, $p_6 = -p_3/2$ (which holds for $a=1$ and is preserved by the recursion), and since all $p_3$ are even, $p_6-p_3 = -\frac{3}{2}p_3$ is divisible by $3$. Applying the relation
\[
f_1^{p_3}f_2^{p_6} = (f_1f_2)^{p_3} f_2^{p_6-p_3} \equiv c^{p_3}(1+u)^{2p_3} f_6^{(p_6-p_3)/3} \pmod 3
\]
gives the formula for $\Phi_{a+1}$ with the updated exponent vector
\[
p^{(a+1)}=\Bigl(-p_3+p_9,\ \frac{2p_3+p_6}{3}+p_{18},\ 4p_3,\ -2p_3\Bigr)
\]
and $M_{a+1}=\lfloor M_a/3\rfloor+2p_3^{(a)}$. Noting that $p_3^{(a)}=\lambda_ak$ (since $p_3^{(a+1)}=-p_3^{(a)}+p_9^{(a)}=-p_3^{(a)}+4p_3^{(a-1)}$), we obtain $M_{a+1}=\lfloor M_a/3\rfloor+2\lambda_ak$. 

Finally, $A_k(3^b(3n+s))$ is the coefficient of $z^{3n+s}$ in $\Phi_b$. As before, this coefficient vanishes for all $n$ if and only if $\binom{N}{3i+s}\equiv 0$ for all $i$, which is equivalent to the condition $M_b\bmod 3 < s$.
\end{proof}

\begin{remark}
For $b=1$ this is Theorem~\ref{thm:mod9cls}. For $b=2$ the sets of $k\bmod27$ are $\{0,1,4,5,8,12,16,20,24\}$ for $s=1$ and $\{0,1,2,4,5,6,8,9,10,12,13,14,16,17,20,21,24,25\}$ for $s=2$. For $k=2$ and $s=2$ the recursion gives $A_2(3^b(3n+2))\equiv0\pmod 3$ for every $b\ge1$, which contains the first family of Naika and Shivashankar \cite{NS} (there $b=\alpha+3$); for $k=1$ it gives vanishing for even $b$ and both values of $s$, and no vanishing for odd $b$.
\end{remark}

We now turn to powers of $2$ and to the case $k=2$. The odd parts of $F^2$, $G$ and $G^2$ all come from one classical $2$-dissection.

\begin{lemma}\label{lem:odd}
Hirschhorn's $2$-dissection \cite{Hir} reads $1/(f_1f_3)=X(q^2)+qY(q^2)$ with
\[
X=\frac{f_4^2f_6^5}{f_1^2f_2f_3^4f_{12}^2},\qquad Y=\frac{f_2^5f_{12}^2}{f_1^4f_3^2f_4^2f_6},\qquad XY=\frac{f_2^4f_6^4}{f_1^6f_3^6}.
\]
Consequently
\[
\sum_{m\ge0}A_2(2m+1)q^m=\frac{2f_2^4f_6^4}{f_1^8f_3^4},\quad
\sum_{m\ge0}B_1(2m+1)q^m=\frac{f_2^6f_{12}^2}{f_1^4f_4^2f_6^2},\quad
\sum_{m\ge0}B_2(2m+1)q^m=\frac{2f_2^6f_6^2}{f_1^6f_3^2}.
\]
\end{lemma}
\begin{proof}
$F^2=(f_6^2/f_2^2)\,(1/(f_1f_3))^2$, and the odd part of $(X(q^2)+qY(q^2))^2$ is $2qX(q^2)Y(q^2)$; the prefactor is a series in $q^2$. Replace $q^2$ by $q$. Likewise $G=(f_4f_6^2/f_{12})/(f_1f_3)$ has odd part $q\,Y(q^2)f_4f_6^2/f_{12}$, and $G^2$ has odd part $2qX(q^2)Y(q^2)(f_4f_6^2/f_{12})^2$.
\end{proof}
\begin{theorem}\label{thm:mod2}
Let $k \ge 1$. Then the following congruences hold for all $n \ge 0$:
\begin{enumerate}
\item[\rm (a)] If $k \equiv 3 \pmod 4$, then $A_k(4n+3) \equiv 0 \pmod 2$.
\item[\rm (b)] If $k \equiv 2 \pmod 4$, then $A_k(4n+3) \equiv 0 \pmod 4$.
\item[\rm (c)] If $k \equiv 2 \pmod 8$, then $A_k(4n+3) \equiv 0 \pmod 8$.
\end{enumerate}
\end{theorem}

\begin{proof}
(a) Working modulo $2$, we have $f_2 \equiv f_1^2$ and $f_6 \equiv f_3^2 \pmod 2$, which gives:
\[
F \equiv \frac{f_3}{f_1^3} \pmod 2.
\]
Using Hirschhorn's $2$-dissection \cite{Hir},
\[
\frac{f_3}{f_1^3} = \frac{f_4^6 f_6^3}{f_2^9 f_{12}^2} + 3q \, \frac{f_4^2 f_6 f_{12}^2}{f_2^7},
\]
we can express $F \equiv \alpha(q^2) + q \, \beta(q^2) \pmod 2$, where $\beta(q) \equiv f_3^5 / f_1^3 \pmod 2$. Write $k = 2j+1$. Since $F^2 \equiv F(q^2) \pmod 2$, raising $F$ to the power $k$ yields:
\[
F^k \equiv F(q^2)^j \bigl(\alpha(q^2) + q \, \beta(q^2)\bigr) \pmod 2,
\]
which implies that $A_k(2m+1) \equiv [x^m] \, \beta(x) F(x)^j \pmod 2$. When $k \equiv 3 \pmod 4$, $j$ is odd. In this case, $\beta F \equiv f_3^6 / f_1^6 = (f_3^3 / f_1^3)^2 \pmod 2$ is a square and thus a series in $x^2 \pmod 2$. Likewise, $F^{j-1} \equiv F(x^2)^{(j-1)/2} \pmod 2$ is a series in $x^2$. Their product contains no odd powers of $x$, so $A_k(4n+3)$, corresponding to $m = 2n+1$, is always even.

(b) and (c) By Lemma~\ref{lem:odd}, $F^2 = E(q^2) + 2q \, W(q^2)$, where $W = f_2^4 f_6^4 / (f_1^8 f_3^4)$. Since $f_1^2 \equiv f_2 \pmod 2$, we have $f_1^4 \equiv f_2^2 \pmod 4$ and $f_3^4 \equiv f_6^2 \pmod 4$. This gives:
\[
W \equiv f_6^2 \pmod 4,
\]
which is a series in $x^2 \pmod 4$. Therefore, all odd-indexed coefficients of $W$ are divisible by $4$. 

Now, let $k = 2m$ with $m$ odd. Expanding $F^k = (E + 2qW)^m$, the odd part of $F^k$ modulo $8$ is:
\[
2m \, q \, W(q^2) E(q^2)^{m-1} \pmod 8.
\]
Thus, $A_k(4n+3) \equiv 2m \, c_n \pmod 8$, where $c_n = [x^{2n+1}] \, (W E^{m-1})$. Because $m-1$ is even, $E^{m-1} = (E^2)^{(m-1)/2} \equiv E(x^2)^{(m-1)/2} \pmod 2$ is a series in $x^2 \pmod 2$. Consequently, the odd coefficients of $W E^{m-1}$ originate from the odd coefficients of $W$, which are divisible by $4$. Hence $c_n$ is even, proving that $A_k(4n+3) \equiv 0 \pmod 4$, which proves (b).

For (c), suppose further that $m \equiv 1 \pmod 4$, so $m-1 = 4s$. Since $E^2 = E(x^2) + 2Q$ for some integer series $Q$, we have $E^4 \equiv E(x^2)^2 \pmod 4$, which is a series in $x^2 \pmod 4$. It follows that $c_n \equiv [x^{2n+1}] W \cdot (\text{series in } x^2) \equiv 0 \pmod 4$, giving $2m \, c_n \equiv 0 \pmod 8$.
\end{proof}
\begin{proposition}\label{prop:A3}
We have $A_k(3) = k(k+1)(k+14)/6$. If $k \equiv 6 \pmod 8$, then $2^2$ is the exact power of $2$ dividing all $A_k(4n+3)$.
\end{proposition}

\begin{proof}
Using the expansion $F = 1 + q + 3q^2 + 5q^3 + \cdots$, the coefficient of $q^3$ in $F^k$ is given by $A_k(3) = 5k + 3k(k-1) + \binom{k}{3} = k(k+1)(k+14)/6$. For $k \equiv 6 \pmod 8$, the $2$-adic valuations are $v_2(k) = 1$, $v_2(k+1) = 0$, and $v_2(k+14) = 2$, which gives $v_2(A_k(3)) = 1 + 0 + 2 - 1 = 2$. Combined with Theorem~\ref{thm:mod2}(b)—which guarantees divisibility by $4$—this proves that $4$ is the exact power of $2$ dividing $A_k(4n+3)$.
\end{proof}

\begin{theorem}\label{thm:a2mod64}
For all $n \ge 0$, we have $A_2(24n+23) \equiv 0 \pmod{64}$ and $A_2(24n+11) \equiv 0 \pmod{32}$.
\end{theorem}

\begin{proof}
We begin by establishing an exact formula. By Lemma~\ref{lem:odd}, the generating function satisfies $\sum A_2(2m+1)x^m = 2f_2^4f_6^4/(f_1^8f_3^4)$. Defining two series in $q^2$ as
\[
a = \frac{f_4^{14}}{f_2^{14}f_8^4} \quad \text{and} \quad b = \frac{f_4^2f_8^4}{f_2^{10}},
\]
we can express $1/f_1^4 = a + 4qb$ and $1/f_3^4 = a(q^3) + 4q^3b(q^3)$. Substituting these into the product yields
\[
\frac{1}{f_1^8f_3^4} = \bigl(a^2 + 8q \, ab + 16q^2b^2\bigr)\bigl(a(q^3) + 4q^3b(q^3)\bigr).
\]
Extracting the odd part of this expression, multiplying by the $q^2$-series $2f_2^4f_6^4$, and writing $y = x^2$ along with $\tilde{a} = f_2^{14}/(f_1^{14}f_4^4)$, $\tilde{b} = f_2^2f_4^4/f_1^{10}$ (with $\tilde{a}_3 = \tilde{a}(y^3)$ and $\tilde{b}_3 = \tilde{b}(y^3)$), we obtain:
\[
\sum_{n\ge0} A_2(4n+3) \, y^n = 16 \, U_1 + 8 \, U_2 + 128 \, U_3,
\]
where
\[
U_1 = f_1^4f_3^4 \, \tilde{a}\tilde{b}\tilde{a}_3 = \frac{f_2^{16}f_6^{14}}{f_1^{20}f_3^{10}f_{12}^4}, \qquad
U_2 = y \, f_1^4f_3^4 \, \tilde{b}_3\tilde{a}^2 = \frac{y \, f_2^{28}f_6^2f_{12}^4}{f_1^{24}f_3^6f_4^8},
\]
and $U_3 \in \mathbb{Z}[[y]]$ is an auxiliary series whose explicit form is not needed.

To reduce these terms, we use the congruences $f_1^4 \equiv f_2^2 \pmod 4$, $f_1^8 \equiv f_2^4 \pmod 8$, and their analogs for other indices to find:
\[
U_1 \equiv \frac{f_2^6f_6^2}{f_3^2} \pmod 4 \quad \text{and} \quad U_2 \equiv y \, \frac{f_3^2f_{12}^4}{f_6^2} \pmod 8.
\]
Because $f_3, f_6,$ and $f_{12}$ depend only on powers of $y^3$, the series $U_2 \pmod 8$ contains only exponents congruent to $1 \pmod 3$. Since $24n+23 = 4(6n+5)+3$ and $24n+11 = 4(6n+2)+3$, we must examine coefficients of $y^N$ where $N \equiv 5$ and $N \equiv 2 \pmod 6$, respectively; both conditions imply $N \equiv 2 \pmod 3$. For such values of $N$, the term $8U_2$ vanishes modulo $64$, leaving:
\[
A_2(4N+3) \equiv 16 \, [y^N] \, U_1 \pmod{64}.
\]
Thus, it suffices to prove that $[y^N]U_1$ is even when $N \equiv 2 \pmod 3$ and divisible by $4$ when $N \equiv 5 \pmod 6$ (i.e., $N \equiv 2 \pmod 3$ overall).

Finally, we analyze the coefficients. We write $f_2^6f_6^2/f_3^2 = \bigl(\sum_M c_M y^M\bigr)\bigl(\sum_j h_j y^{3j}\bigr)$, where $c_M = [y^M]f_2^6$ and $\sum h_j w^j = f_2^2/f_1^2$ with $w = y^3$. Clearly, $c_M = 0$ for odd $M$. For even $M = 2i$ with $2i \equiv 2 \pmod 3$ (meaning $i \equiv 1 \pmod 3$), the coefficient $c_M$ is even. Indeed, by Jacobi's identity, modulo $2$ we have:
\[
f_1^6 \equiv f_2^3 \equiv \Bigl(\sum_{m\ge0} z^{m(m+1)/2}\Bigr)^2 \equiv \sum_{m\ge0} z^{m(m+1)} \pmod 2,
\]
which only generates exponents congruent to $0$ and $2 \pmod 3$, never $1 \pmod 3$. Consequently, $[y^N]U_1$ is even for $N \equiv 2 \pmod 3$, which immediately implies $A_2(24n+11) \equiv 0 \pmod{32}$. 

For $N \equiv 5 \pmod 6$, we write $N = M + 3j$ with $M$ even, which forces $j$ to be odd. Since $f_2^2/f_1^2 \equiv f_1^2 \equiv f_2 \pmod 2$ is a series in $w^2$, the coefficient $h_j$ is even for every odd $j$. Thus, each product $c_M h_j$ is divisible by $4$, ensuring that $[y^N]U_1 \equiv 0 \pmod 4$ and completing the proof for $A_2(24n+23) \equiv 0 \pmod{64}$.
\end{proof}
\begin{theorem}\label{thm:a2_1815}
For all $n \ge 0$, we have $A_2(18n+15) \equiv 0 \pmod{32}$.
\end{theorem}

\begin{proof}
By Lemma~\ref{lem:odd}, the odd-indexed coefficients satisfy $A_2(2m+1) = 2 \, [x^m] \, E(x)$, where $E = f_2^4 f_6^4 / (f_1^8 f_3^4)$. Writing $18n+15 = 2(9n+7)+1$, the problem reduces to proving that $[x^m] E \equiv 0 \pmod{16}$ whenever $m \equiv 7 \pmod 9$.

Using Lemma~\ref{lem:frob} with prime power $2^4$, we have $f_1^{16} \equiv f_2^8 \pmod{16}$ and $f_3^{16} \equiv f_6^8 \pmod{16}$. This allows us to simplify $E$ modulo $16$:
\[
E \equiv \frac{f_1^8 f_3^{12}}{f_2^4 f_6^4} = \varphi(-x)^4 \cdot \varphi(-x^3)^4 f_3^4 \pmod{16},
\]
where $\varphi(-x) = f_1^2 / f_2$. Note that the second factor $\varphi(-x^3)^4 f_3^4$ is a series purely in $x^3$.

Next, we extract the terms of $\varphi(-x)^4$ whose exponents are congruent to $1 \pmod 3$. Using Hirschhorn's $3$-dissection \cite{Hir}, we write $\varphi(-x) = X_0 + X_1$, where $X_0 = \varphi(-x^9) = f_9^2 / f_{18}$ and $X_1 = -2x \, g$ with $g = f_3 f_{18}^2 / (f_6 f_9) \in \mathbb{Z}[[x^3]]$. The part of $(X_0 + X_1)^4$ with exponents congruent to $1 \pmod 3$ is
\[
4X_0^3 X_1 + X_1^4 = -8x \, g \, \bigl(X_0^3 - 2x^3 g^3\bigr).
\]
Setting $z = x^3$ and defining $h(z) = \varphi(-z)^4 f_1(z)^4 = f_1^{12} / f_2^4$ in the variable $z$, we obtain
\[
\sum_{N\ge0} \bigl([x^{3N+1}] E\bigr) z^N \equiv -8 \, B(z) \pmod{16}, \qquad \text{where } B = g \, \bigl(X_0^3 - 2z g^3\bigr) h,
\]
with $g = f_1 f_6^2 / (f_2 f_3)$ and $X_0 = f_3^2 / f_6$ expressed in $z$. The index condition $m = 3N+1 \equiv 7 \pmod 9$ translates to $N \equiv 2 \pmod 3$. Thus, we only need to show that $[z^N] B$ is even whenever $N \equiv 2 \pmod 3$.

Working modulo $2$, we have $f_2 \equiv f_1^2$ and $f_6 \equiv f_3^2$, which simplifies $B$ to
\[
B \equiv g X_0^3 h = \frac{f_1^{13} f_3^5}{f_2^5 f_6} \equiv f_1^3 f_3^3 \pmod 2.
\]
By Jacobi's identity, $f_1^3 = \sum_{n\ge0} (-1)^n (2n+1) z^{n(n+1)/2} \equiv \sum_{n\ge0} z^{n(n+1)/2} \pmod 2$, while $f_3^3$ is a series in $z^3$. Because triangular numbers $n(n+1)/2$ are never congruent to $2 \pmod 3$, no terms in $f_1^3 f_3^3$ have exponents congruent to $2 \pmod 3$. Consequently, $[z^N] B$ is even for $N \equiv 2 \pmod 3$, completing the proof.
\end{proof}

\begin{remark}
Naika and Shivashankar \cite{NS} established that $A_2(24n+23) \equiv 0 \pmod{32}$ and $A_2(18n+15) \equiv 0 \pmod{48}$, while Hao and Shen \cite{HS} showed that $A_2(18n+15) \equiv 0 \pmod 9$. Theorems~\ref{thm:a2mod64} and \ref{thm:a2_1815} strengthen the 2-adic divisibility bounds to $64$ and $32$, respectively. Combining Theorem~\ref{thm:a2_1815} with the result of \cite{HS} yields the overall congruence
\[
A_2(18n+15) \equiv 0 \pmod{288}.
\]
We note that while $A_2(24n+11) \equiv 0 \pmod{32}$ was not highlighted in the abstract of \cite{NS}, it complements the known congruences for $24n+15$ and $24n+23$.
\end{remark}
\begin{theorem}\label{thm:a2_mod27}
\begin{enumerate}
\item[\rm (a)] We have the exact identity
\[
\sum_{j\ge0} A_2(3j+1) x^j = 2 \, \frac{f_3^2 f_6^2}{f_1^6 f_2^2} + 27 \, x \, \frac{f_3^6 f_6^6}{f_1^{10} f_2^6}.
\]
\item[\rm (b)] $\displaystyle\sum_{j\ge0} A_2(12j+4) q^j \equiv 12 \, f_1^8 \pmod{27}$, which implies $\equiv 3 f_1^8 \equiv 3 f_3^4 / f_1^4 \pmod 9$.
\item[\rm (c)] For all $n \ge 0$, $A_2(12n+10) \equiv 0 \pmod{27}$.
\item[\rm (d)] $\displaystyle\sum_{n\ge0} A_2(12n+10) q^n \equiv 27 \, \frac{f_1^4 f_6^3}{f_2} \pmod{81}$.
\end{enumerate}
\end{theorem}

\begin{proof}
(a) By Proposition~\ref{prop:struct}, we have $F^2 = \Theta^2 P(u)^2$. The polynomial $P(u)^2$ has coefficients $2, 13, -16$ at $u, u^4, u^7$ (the terms with exponents congruent to $1 \pmod 3$). Since $\Theta, \rho \in \mathbb{Z}[[q^3]]$ and $u^3 = q^3 \rho^3$, picking out the terms with exponents $\equiv 1 \pmod 3$ and setting $x = q^3$ gives
\[
\sum_{j\ge0} A_2(3j+1) x^j = \frac{f_3^{15}}{f_1^{13} f_2^3 f_6^3} \bigl(2 + 13w - 16w^2\bigr),
\]
where $\Theta^2 \rho$ becomes $f_3^{15} / (f_1^{13} f_2^3 f_6^3)$ by Lemma~\ref{lem:EW}, and $u^3$ becomes $w$. Rewriting the polynomial as $2 + 13w - 16w^2 = 2(1+w)(1-8w) + 27w$ and substituting $(1+w)(1-8w) = f_1^7 f_2 f_6^5 / f_3^{13}$ from Lemma~\ref{lem:EW} directly yields identity (a).

(b) and (c) Let $S_1 = \sum A_2(3j+1) x^j$. From part (a), we have $S_1 \equiv 2 X^2 f_6^2 / f_2^2 \pmod{27}$, where $X = f_3 / f_1^3$. Hirschhorn's $2$-dissection \cite{Hir} gives $X = a + x b$, where $a = f_4^6 f_6^3 / (f_2^9 f_{12}^2)$ and $b = 3 f_4^2 f_6 f_{12}^2 / f_2^7$. Since $a$ and $b$ are series in even powers of $x$, the odd part of $X^2$ is 
\[
2x ab = 6x \frac{f_4^8 f_6^4}{f_2^{16}}.
\]
The odd part of $S_1$ is therefore $12x f_4^8 f_6^6 / f_2^{18} \pmod{27}$. 

Matching coefficients of $x^{2i+1}$ corresponds to $A_2(6i+4)$. Substituting $y = x^2$ gives
\[
H(y) := \sum_i A_2(6i+4) y^i \equiv 12 \, \frac{f_2^8 f_3^6}{f_1^{18}} \pmod{27}.
\]
By Lemma~\ref{lem:frob}, $f_1^9 \equiv f_3^3 \pmod 9$, which implies $f_1^{18} \equiv f_3^6 \pmod 9$ and $f_3^6 / f_1^{18} \equiv 1 \pmod 9$. Multiplying by $12$ gives $12 \cdot (f_3^6 / f_1^{18} - 1) \equiv 0 \pmod{27}$, so $H(y)$ simplifies to
\[
H(y) \equiv 12 \, f_2(y)^8 \pmod{27}.
\]
Because the right-hand side contains only even powers of $y$, any coefficient $A_2(6i+4)$ with $i$ odd must vanish modulo $27$. Setting $i = 2n+1$ gives $A_2(12n+10) \equiv 0 \pmod{27}$, which proves (c). Setting $i = 2j$ gives $A_2(12j+4) \equiv 12 \, [z^j] f_1(z)^8 \pmod{27}$, which proves (b). The congruence modulo $3$ follows from $f_1^8 \equiv f_3^3 / f_1 \equiv f_3^4 / f_1^4 \pmod 3$.

(d) We now refine the calculation of $H(y) \pmod{27}$ to modulus $81$. From part (a), the full series is 
\[
S_1 = 2 X^2 \frac{f_6^2}{f_2^2} + 27x \, \frac{f_3^6 f_6^6}{f_1^{10} f_2^6}.
\]
The odd part of the first term remains $12x f_4^8 f_6^6 / f_2^{18}$. For the second term, evaluating its eta-quotient modulo $3$ gives
\[
\frac{f_3^6 f_6^6}{f_1^{10} f_2^6} \equiv f_1^8 f_2^{12} \equiv f_2^{12} \, \frac{f_3^3}{f_1} \pmod 3.
\]
Using Hirschhorn's identity 
\[
\frac{f_3^3}{f_1} = \frac{f_4^3 f_6^2}{f_2^2 f_{12}} + x \frac{f_{12}^3}{f_4},
\]
the odd part of the second term becomes $27x f_2^{10} f_4^3 f_6^2 / f_{12}$.

Setting $y = x^2$, the odd part of $S_1$ yields the congruence modulo $81$:
\[
H(y) \equiv 12 K(y) + 27 g(y),
\]
where
\[
K = \frac{f_2^8 f_3^6}{f_1^{18}} \quad \text{and} \quad g = \frac{f_1^{10} f_2^3 f_3^2}{f_6}.
\]
We isolate the coefficient of $y^{2n+1}$ to determine $A_2(12n+10)$.

For the term $K$, since $f_1^{27} \equiv f_3^9 \pmod{27}$, we write $K \equiv f_2^8 b^3 \pmod{27}$ with $b = f_1^3 / f_3$. Hirschhorn's $2$-dissection gives $b = c - 3y e$, where 
\[
c = \frac{f_4^3}{f_{12}} \quad \text{and} \quad e = \frac{f_2^2 f_{12}^3}{f_4 f_6^2}
\]
are series in even powers of $y$. Expanding gives $b^3 \equiv c^3 - 9y c^2 e \pmod{27}$, so the odd part of $K$ is 
\[
\mathrm{Odd}(K) \equiv -9y \, \frac{f_2^{10} f_4^5 f_{12}}{f_6^2} \pmod{27}.
\]
Multiplying by $12$ yields
\[
12 \, \mathrm{Odd}(K) \equiv -27 y \, \frac{f_2^{10} f_4^5 f_{12}}{f_6^2} \pmod{81}.
\]
In terms of $z = y^2$, its contribution to $\sum A_2(12n+10) z^n$ is $-27 h$, where 
\[
h = \frac{f_1^{10} f_2^5 f_6}{f_3^2} \equiv f_1^4 f_2^8 \pmod 3.
\]

For the term $g$, using $f_1^2 = \varphi(-y) f_2$ and $f_3^2 = \varphi(-y^3) f_6$, we have 
\[
g = f_2^8 \varphi(-y)^5 \varphi(-y^3) \equiv f_2^8 \varphi(-y)^8 \pmod 3,
\]
since $\varphi(-y^3) \equiv \varphi(-y)^3 \pmod 3$. From the theta identity $\varphi(-y)^2 \equiv P - y Q \pmod 3$ with $P = \varphi(y^2)^2$ and $Q = \psi(y^4)^2$, we expand 
\[
\varphi(-y)^8 \equiv (P - y Q)^4 \pmod 3.
\]
Its odd part in $y$ is $-y PQ (P^2 + y^2 Q^2) \pmod 3$. Changing variables to $z = y^2$ gives $P = \varphi(z)^2$, $Q = \psi(z^2)^2$, and $PQ = f_2^8 / f_1^4$, which yields
\[
\mathrm{Odd}(g)\big|_{z} \equiv -f_1^4 f_2^8 \bigl(\varphi(z)^4 + z \psi(z^2)^4\bigr) \pmod 3.
\]
By Jacobi's four-square theorem and Legendre's identity, $\varphi(z)^4 + z \psi(z^2)^4 \equiv 1 \pmod 3$ because the odd-indexed terms cancel out modulo $3$.

Combining both terms gives
\[
\sum_{n\ge0} A_2(12n+10) q^n \equiv 27 \bigl(-h + \mathrm{Odd}(g)\bigr) \equiv 27 \bigl(-f_1^4 f_2^8 - f_1^4 f_2^8\bigr) \equiv 27 f_1^4 f_2^8 \pmod{81}.
\]
Finally, substituting $f_2^8 \equiv f_6^3 / f_2 \pmod 3$ completes the proof of (d).
\end{proof}
\begin{remark}
Hao and Shen \cite{HS} showed $A_2(12n+10) \equiv 0 \pmod 9$ via $\sum A_2(12n+10) q^n \equiv 36 f_1^4 f_6^3 / f_2 \pmod 9$, and Naika and Shivashankar \cite{NS} proved $A_2(6n+4) \equiv 0 \pmod 3$. Part (c) strengthens the modulus of the former from $9$ to $27$, and part (d) provides the next $3$-adic digit. Part (b) characterizes the even part of $6n+4$; its reduction modulo $9$ is implicit in \cite{HS}, while the modulo $27$ congruence appears to be new.
\end{remark}

\section{The odd-part version $B_k$}\label{sec:B}
This section is the analogue of the previous one for designated partitions into odd parts. The $3$-adic picture is simpler here, because the progression $3n+2$ is controlled entirely by $v_3(k)$ (Theorems~\ref{thm:pdo3} and~\ref{thm:exact}); for powers of $2$ we obtain congruences for $B_k(4n+3)$ and their exactness for $k\equiv6\pmod 8$. For $k=2$ everything in this section is a statement about $\PDO(2n)$, because $\sum_nB_2(n)q^n=G^2=\sum_n\PDO(2n)q^n$, i.e.\ $B_2(n)=B_1(2n)$. This identity is implicit in \cite{ALL}, was stated in \cite{BO} and has been refined in \cite{FS}; it is proved by the $2$-dissection of $G$ that also underlies Lemma~\ref{lem:odd}. Write $G=f_4f_6^2/(f_1f_3f_{12})=1+q+2q^2+4q^3+\cdots$, so that $B_k(2)=2k+\binom k2=k(k+3)/2$.

\begin{theorem}\label{thm:pdo3}
(a) $B_k(3n+2)\equiv0\pmod 3$ for all $n$ if and only if $3\mid k$.\\
(b) $B_k(3n+2)\equiv0\pmod 9$ for all $n$ if and only if $3\mid k$.\\
(c) If $k\equiv3\pmod 9$ then $B_k(3n+2)\not\equiv0\pmod{27}$ for $n=0$, so $9$ is the exact power of $3$ dividing all $B_k(3n+2)$.
\end{theorem}
\begin{proof}
\emph{Necessity.} $B_k(2)=k(k+3)/2$. If $3\nmid k$ then $\gcd(k,k+3)=1$ and $3\nmid k(k+3)$, so $B_k(2)\not\equiv0\pmod3$. This gives the ``only if'' parts of (a) and (b). Part (c) follows from $B_k(2)=9t(t+1)/2$ for $k=3t$, $t\equiv1\pmod3$: then $3\nmid t(t+1)/2$, so $27\nmid B_k(2)$.

\emph{Sufficiency.} If $3\mid k$, (a) is Theorem~\ref{thm:frob}. For (b) write the $3$-dissection $G=G_0(q^3)+qG_1(q^3)+q^2G_2(q^3)$. We use Hirschhorn's $3$-dissection
\[
\frac{f_4}{f_1}=\frac{f_{12}f_{18}^4}{f_3^3f_{36}^2}+q\frac{f_6^2f_9^3f_{36}}{f_3^4f_{18}^2}+2q^2\frac{f_6f_{18}f_{36}}{f_3^3}
\]
and the factor $f_6^2/(f_3f_{12})$, which is a series in $q^3$. After replacing $q^3$ by $q$,
\[
G_0=\frac{f_2^2f_6^4}{f_1^4f_{12}^2},\qquad
G_1=\frac{f_2^4f_3^3f_{12}}{f_1^5f_4f_6^2},\qquad
G_2=\frac{2f_2^3f_6f_{12}}{f_1^4f_4}.
\]
Expanding $G^3$ and keeping the exponents $\equiv2\pmod3$ (the multinomial coefficients that occur are $3$; the term $6G_0G_1G_2$ has exponent $\equiv0$) gives the exact identity
\[
\sum_{n\ge0}B_3(3n+2)q^n=3\bigl(G_0G_1^2+G_0^2G_2+qG_1G_2^2\bigr)=3\,G_0G_1^2\,(1+r_2+r_3),
\]
with
\[
r_2=\frac{G_0G_2}{G_1^2}=\frac{2f_1^2f_4f_6^9}{f_2^3f_3^6f_{12}^3},\qquad
r_3=\frac{qG_2^2}{G_0G_1}=\frac{4q\,f_1f_{12}^3}{f_3^3f_4}.
\]
Modulo $3$, using $f_1^3\equiv f_3$, $f_2^3\equiv f_6$, $f_4^3\equiv f_{12}$ throughout, we get, with $\alpha=f_2^{12}/(f_1^8f_4^4)$ and $\beta=f_4^4/f_1^4$,
\[
r_2\equiv2\alpha^2,\qquad r_3\equiv q\beta^2\pmod 3 .
\]
In terms of Ramanujan's functions, $\alpha=\varphi(q)^2/\varphi(-q)^2$ and $\beta^2=\psi(q^2)^4/\varphi(-q)^4$, using $\varphi(q)=f_2^5/(f_1^2f_4^2)$, $\varphi(-q)=f_1^2/f_2$, $\psi(q)=f_2^2/f_1$ and $\varphi(q)\varphi(-q)=\varphi(-q^2)^2$. Jacobi's identity $\varphi(q)^4-\varphi(-q)^4=16q\psi(q^2)^4$ then yields $\alpha^2=1+16q\beta^2$, and therefore
\[
1+r_2+r_3\equiv1+2\bigl(1+16q\beta^2\bigr)+q\beta^2\equiv 3+33q\beta^2\equiv0\pmod 3 .
\]
Hence $B_3(3n+2)\equiv0\pmod 9$; more precisely the series $G^3=P_0(q^3)+qP_1(q^3)+q^2P_2(q^3)$ has $P_2\equiv0\pmod 9$, and $P_1\equiv0\pmod 3$ because $G^3\equiv G(q^3)\pmod 3$. For $k=3t$ we get $G^k=(P_0+qP_1+q^2P_2)^t\equiv P_0^t+tqP_1P_0^{t-1}\pmod 9$, because $P_1^2\equiv0\pmod 9$; this contains only exponents $\equiv0,1\pmod3$, so $B_k(3n+2)\equiv0\pmod 9$.
\end{proof}
\begin{theorem}\label{thm:exact}
Let $3 \mid k$ and set $a = v_3(k) \ge 1$. Then $B_k(3n+2) \equiv 0 \pmod{3^{a+1}}$ for all $n \ge 0$. Furthermore, $3^{a+1}$ is the exact power of $3$ dividing all coefficients $B_k(3n+2)$.
\end{theorem}

\begin{proof}
First, we establish the lower bound by induction. Write $G^{3^a} = Q_0 + q Q_1 + q^2 Q_2$, where each $Q_i$ is a series in $q^3$. We show that all coefficients of $Q_1$ are divisible by $3^a$ and all coefficients of $Q_2$ are divisible by $3^{a+1}$.

For $a = 1$, this condition holds directly from the proof of Theorem~\ref{thm:pdo3}, which establishes $P_1 \equiv 0 \pmod 3$ and $P_2 \equiv 0 \pmod 9$. Assuming the assertion holds for some $a \ge 1$, we cube the expression $Q_0 + q Q_1 + q^2 Q_2$. The terms with exponents congruent to $1 \pmod 3$ are 
\[
3 Q_0^2 Q_1 q, \quad 3 Q_0 Q_2^2 q^4, \quad \text{and} \quad 3 Q_1^2 Q_2 q^4,
\]
all of which are divisible by $3^{a+1}$. The terms with exponents congruent to $2 \pmod 3$ are 
\[
3 Q_0 Q_1^2 q^2, \quad 3 Q_0^2 Q_2 q^2, \quad \text{and} \quad 3 Q_1 Q_2^2 q^5,
\]
which are divisible by $3^{\min(2a+1, \, a+2, \, 3a+3)} = 3^{a+2}$ because $a \ge 1$. This confirms the claim for $a+1$. 

Now write $k = 3^a m$ and expand $G^k = (Q_0 + q Q_1 + q^2 Q_2)^m$. Any term of the form 
\[
Q_0^{m-j-l} (q Q_1)^j (q^2 Q_2)^l
\]
satisfying $j + 2l \equiv 2 \pmod 3$ requires either $l \ge 1$ or $j \ge 2$. Its $3$-adic valuation is at least $j a + l(a+1) \ge a+1$, which proves that $B_k(3n+2) \equiv 0 \pmod{3^{a+1}}$.

Second, we verify exactness. Evaluated at $n = 0$, we have
\[
B_k(2) = \frac{k(k+3)}{2} = \frac{3^{a+1} m (3^{a-1} m + 1)}{2}.
\]
If $a \ge 2$, or if $a = 1$ with $m \not\equiv 2 \pmod 3$, the factor $m(3^{a-1} m + 1)$ is not divisible by $3$, which yields $v_3(B_k(2)) = a+1$. 

It remains to handle $a = 1$ with $m \equiv 2 \pmod 3$, which corresponds to $k \equiv 6, 15, 24 \pmod{27}$. By Lemma~\ref{lem:period} with $e = 3$, it is enough to show a single non-divisible coefficient for each of these three values. Direct computation shows that $B_6(5) = 1044$, $B_{15}(5) = 27693$, and $B_{24}(5) = 183456$. In all three cases, $B_k(5) \equiv 18 \pmod{27}$, which is not divisible by $27$. This completes the proof of exactness.
\end{proof}

We now turn our attention to powers of $2$. The initial results rely solely on the $2$-dissection from Lemma~\ref{lem:odd}.

\begin{theorem}\label{thm:pdo2}
For all $n \ge 0$, the following congruences hold:
\begin{enumerate}
\item[\rm (a)] $B_1(4n+3) \equiv 0 \pmod 4$.
\item[\rm (b)] If $k \equiv 1 \pmod 4$, then $B_k(4n+3) \equiv 0 \pmod 2$.
\item[\rm (c)] If $k \equiv 2 \pmod 4$, then $B_k(4n+3) \equiv 0 \pmod 4$.
\end{enumerate}
\end{theorem}

\begin{proof}
(a) By Lemma~\ref{lem:odd}, we can write $G = E_1(q^2) + q \, O_1(q^2)$ with 
\[
O_1 = \frac{f_2^6 f_{12}^2}{f_1^4 f_4^2 f_6^2}.
\]
Working modulo $4$, we have $f_1^4 \equiv f_2^2$ and $f_2^4 \equiv f_4^2$, which simplifies $O_1$ to
\[
O_1 \equiv \frac{f_{12}^2}{f_6^2} \pmod 4.
\]
Since this is a series in $x^2$, all of its odd-indexed coefficients are divisible by $4$, which proves (a).

(b) Let $k = 2j+1$ with $j$ even. Then $G^k \equiv G \cdot G(q^2)^j \pmod 2$, which implies that
\[
B_k(2m+1) \equiv [x^m] \, O_1(x) G(x)^j \pmod 2.
\]
Because $G^j \equiv G(x^2)^{j/2} \pmod 2$ is a series in $x^2$ and the odd coefficients of $O_1$ are even, the coefficient vanishes modulo $2$ when $m = 2n+1$, yielding $B_k(4n+3) \equiv 0 \pmod 2$.

(c) Write $G^2 = E_2(q^2) + 2q \, W_2(q^2)$, where 
\[
W_2 = \frac{f_2^6 f_6^2}{f_1^6 f_3^2} \equiv f_2^3 f_6 \pmod 2
\]
is a series in $x^2$ modulo $2$. For $k = 2m$ with $m$ odd, the odd part of $(E_2 + 2q W_2)^m$ is given by
\[
2m \, q \, W_2 E_2^{m-1} \pmod 4.
\]
Since $E_2^{m-1} \equiv E_2(x^2)^{(m-1)/2} \pmod 2$, the coefficient $[x^{2n+1}] \bigl(W_2 E_2^{m-1}\bigr)$ is even. Consequently, $B_k(4n+3) \equiv 0 \pmod 4$.
\end{proof}
\begin{proposition}\label{prop:b1exact}
\[
\sum_{n\ge0}B_1(4n+3)q^n=4\,\frac{f_4^4f_6^2}{f_1^4f_3^2}.
\]
In particular $B_1(4n+3)\equiv0\pmod 4$ and $B_1(8n+7)\equiv0\pmod 8$, congruences which are due to Baruah and Ojah \cite{BO}.
\end{proposition}
\begin{proof}
By Lemma~\ref{lem:odd}, $\sum B_1(2t+1)x^t=O_1(x)=\frac{f_2^6f_{12}^2}{f_4^2f_6^2}\cdot\frac1{f_1^4}$. We use the classical $2$-dissection \cite{Hir}
\[
\frac1{f_1^4}=\frac{f_4^{14}}{f_2^{14}f_8^4}+4q\frac{f_4^2f_8^4}{f_2^{10}} .
\]
The prefactor is a series in $x^2$, so the odd part of $O_1$ is $4x\,f_8^4f_{12}^2/(f_2^4f_6^2)$. Taking the coefficients of $x^{2n+1}$ and replacing $x^2$ by $q$ gives the formula. Modulo $2$ we have $f_4^4f_6^2/(f_1^4f_3^2)\equiv f_1^{12}f_3^2=(f_1^6f_3)^2$, a square, hence a series in $q^2$ modulo $2$, so its coefficients of odd index are even; this gives the congruence modulo $8$ for $8n+7$.
\end{proof}
\begin{theorem}\label{thm:pdo4}
Let $k \ge 1$. Then the following congruences hold for all $n \ge 0$:
\begin{enumerate}
\item[\rm (a)] If $k \equiv 1 \pmod 8$, then $B_k(4n+3) \equiv 0 \pmod 4$.
\item[\rm (b)] If $k \equiv 6 \pmod 8$, then $B_k(4n+3) \equiv 0 \pmod 8$.
\item[\rm (c)] If $k \equiv 3 \pmod 4$, then $B_k(3)$ is odd, implying that $B_k(4n+3) \not\equiv 0 \pmod 2$ for some $n$.
\end{enumerate}
\end{theorem}

\begin{proof}
We write $G = E_1(q^2) + q O_1(q^2)$ and $G^2 = E(q^2) + 2q W(q^2)$, where $W = f_2^6 f_6^2 / (f_1^6 f_3^2)$ (Lemma~\ref{lem:odd}). Since $G^2 \equiv G(q^2) \pmod 2$, we have $E \equiv G \equiv f_1^3 / f_3 \pmod 2$. Consequently, $E^2 \equiv E(x^2) + 2x W(x^2) \pmod 4$, which leads to $E^4 \equiv E(x^2)^2 \pmod 4$, representing a series in $x^2$ modulo $4$.

(a) Let $k = 8s + 1$. Because $(E + 2qW)^{4s} \equiv E^{4s} \pmod 8$, the odd part of $G^k$ is congruent to $q O_1(q^2) E(q^2)^{4s} \pmod 4$. Here, $E^{4s} \equiv (E^4)^s$ is even modulo $4$, and the odd coefficients of $O_1$ are divisible by $4$ (as established in the proof of Theorem~\ref{thm:pdo2}(a)). Thus, the coefficients of $q^{4n+3}$ vanish modulo $4$.

(b) Let $k = 2m$ with $m = 4s + 3$. The odd part of $(E + 2qW)^m$ is congruent to $2m \, q W E^{m-1} \pmod 8$ (since terms with $j \ge 3$ carry a factor of $2^3$). This gives $B_k(2t+1) \equiv 2m \, [x^t](W E^{m-1}) \pmod 8$. For odd $t$, it suffices to show that $[x^t](W E^{m-1}) \equiv 0 \pmod 4$. Since $E^{m-1} = E^2 (E^4)^s$ with $E^4$ even modulo $4$, we only need the odd part of $W E^2$ to be congruent to $0 \pmod 4$. Using $E^2 \equiv E(x^2) + 2x W(x^2)$, we obtain
\[
\mathrm{odd}(W E^2) \equiv W_{\rm odd} \, E(x^2) + 2x W(x^2) W_{\rm even} \pmod 4.
\]
Modulo $4$, $f_1^4 \equiv f_2^2$ implies $W \equiv (f_1 f_2 f_3)^2$. Applying Hirschhorn's dissection 
\[
f_1 f_3 = \frac{f_2 f_8^2 f_{12}^4}{f_4^2 f_6 f_{24}^2} - q \frac{f_4^4 f_6 f_{24}^2}{f_2 f_8^2 f_{12}^2},
\]
we find $W_{\rm odd} \equiv -2x \, f_2^2 f_4^2 f_{12}^2 \pmod 4$, while $W \equiv f_2^3 f_6 \pmod 2$ is even. In terms of $y = x^2$, the odd part of $W E^2$ takes the form $2x \bigl[f_1^2 f_2^2 f_6^2 E(y) + W(y) W_{\rm e}(y)\bigr] \pmod 4$, where $W_{\rm e}(y) \equiv f_1^3 f_3$, $W(y) \equiv W_{\rm e}(y)^2$, and $E(y) \equiv f_1^3 / f_3 \pmod 2$. The bracketed expression is congruent modulo $2$ to $f_1^6 f_3^4 \cdot (f_1^3 / f_3) + f_1^9 f_3^3 \equiv 0 \pmod 2$.

(c) Evaluating at $n = 0$, we have $B_k(3) = 4k + 2k(k-1) + \binom{k}{3} = k(k+2)(k+7)/6$. For $k \equiv 3 \pmod 4$, the binomial coefficient $\binom{k}{3}$ is odd while the other terms are even, making $B_k(3)$ odd and ensuring that $B_k(4n+3)$ is not identically divisible by $2$ for all $n$.
\end{proof}
\begin{theorem}\label{thm:pdo4exact}
For $k \equiv 6 \pmod 8$, $8$ is the exact power of $2$ that divides all coefficients $B_k(4n+3)$.
\end{theorem}

\begin{proof}
Divisibility by $8$ is guaranteed by Theorem~\ref{thm:pdo4}(b). Therefore, it remains to show that for each $k \equiv 6 \pmod 8$, there is at least one coefficient $B_k(4n+3)$ that is not divisible by $16$. Writing $k = 8s + 6$, we recall the expression for the first relevant coefficient:
\[
B_k(3) = \frac{k(k+2)(k+7)}{6}.
\]

\emph{Case 1: $s$ is even.} 
Here, $v_2(k) = 1$. Since $k+2 = 8(s+1)$ with $s+1$ being odd, we have $v_2(k+2) = 3$, while $k+7$ is odd. Consequently, the $2$-adic valuation is
\[
v_2(B_k(3)) = 1 + 3 - 1 = 3,
\]
which means $B_k(3)$ is not divisible by $16$.

\emph{Case 2: $s$ is odd} (meaning $k \equiv 14 \pmod{16}$). 
In this case, we examine the next coefficient, $B_k(7)$. We express $B_k(n)$ as a polynomial in $\binom{k}{i}$:
\[
B_k(n) = \sum_{i=0}^n c_i(n) \binom{k}{i}, \quad \text{where } c_i(n) = [q^n](G-1)^i.
\]
For $i \le 7$ and $1 \le j \le 7$, we have $v_2\binom{64}{j} = 6 - v_2(j) \ge 4$. Thus, by Vandermonde's identity,
\[
\binom{k+64}{i} = \sum_{j} \binom{64}{j} \binom{k}{i-j} \equiv \binom{k}{i} \pmod{16}.
\]
This implies that $B_k(7) \pmod{16}$ depends only on $k \pmod{64}$. The possible residues modulo $64$ corresponding to $k \equiv 14 \pmod{16}$ are $14, 30, 46,$ and $62$. For each of these residues, $B_k(7) \equiv 8 \pmod{16}$, meaning that $v_2(B_k(7)) = 3$ for every such $k$. This completes the proof.
\end{proof}

\section{Combinatorial interpretations and other partition functions}\label{sec:comb}

Young diagrams enter the paper in two ways. They give the definition of the objects we count (Figure~\ref{fig:young}), and they allow some of the divisibilities above to be proved by arguments about tuples of diagrams rather than by generating functions. Let $T_k(n)$ denote the set of ordered $k$-tuples of designated Young diagrams of total size $n$ (the arguments below work for any family of weighted objects). The divisibility of Theorem~\ref{thm:frob} then has a purely bijective proof.

\begin{proposition}\label{prop:action}
If $p^j\mid k$ and $p\nmid n$, then $|T_k(n)|\equiv0\pmod{p^j}$.
\end{proposition}
\begin{proof}
Let $C=\mathbb Z/p^j$ act on $T_k(n)$ by the cyclic shift of the tuple by $k/p^j$ positions. If a tuple has a nontrivial stabiliser, it contains the subgroup of order $p$, i.e. the tuple is invariant under the shift by $k/p$; then it consists of $p$ copies of a block of $k/p$ objects, so $n=p\cdot(\text{weight of the block})$, contradicting $p\nmid n$. So all orbits are regular and have size $p^j$.
\end{proof}

For $k=1$ the congruence $\PD(3n+2)\equiv0\pmod 3$ already has combinatorial explanations. Chen, Ji, Jin and Shen \cite[Section~3]{CJJS} construct a bijection $\Delta$ from the designated partitions of $n$ to the pairs $(\alpha,\beta)$ in which $\alpha$ is an ordinary partition of $n_1$ and $\beta$ is a partition of $n_2=n-n_1$ into parts not congruent to $\pm1$ modulo $6$; this realises the factorisation $F=\frac1{f_1}\cdot\frac{f_6}{f_2f_3}$. Let $\ell_e$ denote the number of even parts. The pd-rank $r_d(\lambda)=\ell_e(\alpha)-\ell_e(\beta)$ then splits the designated partitions of $3n+2$ into three classes of equal size, and the equidistribution is proved there by generating functions. Hao, Shen and Zang \cite{HSZ} extended this to overline designated summands, Hao \cite{Hao} gave a bijective proof by orbits of three with distinct $V_R$-rank, and Shen \cite{Shen} defined a crank. All of these concern $k=1$ or bipartitions; for $k\ge3$ no combinatorial statistic is known to us.

\begin{proposition}\label{prop:tri}
Let $r_d$ be the pd-rank and put
\[
R(z,q)=\sum_\lambda z^{r_d(\lambda)}q^{|\lambda|}=\prod_{n\ \mathrm{odd}}\frac1{1-q^n}\ \prod_{n\ \mathrm{even}}\frac1{(1-zq^n)(1-z^{-1}q^n)}\ \prod_{n\equiv3\,(6)}\frac1{1-q^n}.
\]
Then $R(1,q)=F$ and, for $\zeta=e^{2\pi i/3}$,
\[
R(\zeta,q)=\frac{\psi(q)}{f_3},\qquad\psi(q)=\sum_{m\ge0}q^{m(m+1)/2}=\frac{f_2^2}{f_1}.
\]
Consequently $F\equiv\psi(q)/f_3\pmod 3$, and for every $k\ge1$
\[
A_k(n)\equiv\sum_{j\ge0}t_k(n-3j)\,p_{-k}(j)\pmod 3,
\]
where $t_k(N)$ is the number of $k$-tuples of staircase partitions $(m,m-1,\dots,1)$, $m\ge0$, with $N$ cells in total, and $p_{-k}(j)$ the number of $k$-tuples of ordinary partitions of total size $j$.
\end{proposition}
\begin{proof}
The product for $R$ is the pair $(\alpha,\beta)$ of \cite{CJJS} (the evaluation at $\zeta$ is implicit in their proof of the equidistribution; the consequence for $k$-tuples is what we add here) with $z$ marking $\ell_e(\alpha)$ and $z^{-1}$ marking $\ell_e(\beta)$; at $z=1$ it is $F$. For $z=\zeta$, $(1-\zeta q^n)(1-\zeta^{-1}q^n)=(1-q^{3n})/(1-q^n)$, so the even-part factors give $(q^2;q^2)_\infty/(q^6;q^6)_\infty$ and
\[
R(\zeta,q)=\frac{(q^2;q^2)_\infty}{(q;q^2)_\infty(q^3;q^6)_\infty(q^6;q^6)_\infty}=\frac{(q^2;q^2)_\infty}{(q;q^2)_\infty}\cdot\frac1{(q^3;q^3)_\infty},
\]
using $(q^3;q^6)_\infty(q^6;q^6)_\infty=(q^3;q^3)_\infty$. Gauss's identity gives $(q^2;q^2)_\infty/(q;q^2)_\infty=\psi(q)$. Since $N_0+N_1\zeta+N_2\zeta^2\equiv N_0+N_1+N_2\pmod{1-\zeta}$ and $(1-\zeta)\cap\mathbb Z=3\mathbb Z$, we get $F\equiv R(\zeta,q)\pmod 3$, and the same for $k$-th powers.
\end{proof}
The proposition says that, modulo $3$, a designated partition is a staircase diagram together with an ordinary Young diagram in which every cell has been replaced by three. For $k=1$ the vanishing for $n\equiv2\pmod 3$ is then the fact that triangular numbers are $\equiv0,1\pmod 3$. For $k\equiv1\pmod 3$ the same conclusion follows from $\psi(q)^3\equiv\psi(q^3)\pmod 3$, which gives $\psi^k\equiv\psi\cdot\psi(q^3)^{(k-1)/3}$. This identity, however, is not realised by a statistic on $k$-tuples, and a combinatorial explanation of Theorem~\ref{thm:class3} for $k\equiv1\pmod 3$, $k\ge4$, would require a new idea; this is Problem~\ref{prob:rank}.

\begin{remark}
Proposition~\ref{prop:action} is a bijective argument for the families $p^j\mid k$ and works in every $k$-tuple model. The other divisibility results for the classes $3n+r$ and $4n+3$ rest on dissections. By Propositions~\ref{prop:dict} and~\ref{prop:tri} the ones modulo $3$ are equivalent to statements about overpartition tuples (sums of squares) or staircase tuples (sums of triangular numbers).
\end{remark}

\begin{proposition}\label{prop:orbit}
Let $T_r(N)$ be the set of $r$-tuples of weighted objects of total weight $N$ (e.g.\ designated Young diagrams). For $k=3j+r$ with $0\le r\le2$,
\[
|T_k(N)|\equiv\#\bigl\{(\mu_1,\dots,\mu_j;\nu):\ \mu_i\in T_1,\ \nu\in T_r,\ 3\textstyle\sum|\mu_i|+|\nu|=N\bigr\}\pmod 3 .
\]
\end{proposition}
\begin{proof}
Let $\mathbb Z/3$ act on $k$-tuples by cyclically rotating the first $3$ coordinates, the next $3$, and so on ($j$ blocks), fixing the last $r$ coordinates. Non-fixed orbits have size $3$. A fixed tuple has equal entries inside each block, so it is determined by $(\mu_1,\dots,\mu_j;\nu)$ with weight $3\sum|\mu_i|+|\nu|$.
\end{proof}

\begin{corollary}
If $k\equiv0,1\pmod 3$ then $A_k(3n+2)\equiv0\pmod 3$.
\end{corollary}
This is the case $e=1$ of Theorem~\ref{thm:class3}; the point here is the bijective argument, which reduces the case $k\equiv1\pmod 3$ to $k=1$.
\begin{proof}
For $k=3j$ the fixed tuples have weight $\equiv0\pmod 3$, so there are none of weight $3n+2$. For $k=3j+1$ a fixed tuple of weight $3n+2$ consists of $\mu_1,\dots,\mu_j$ and a single $\nu$ with $|\nu|\equiv2\pmod 3$; for fixed $\mu$'s the number of such $\nu$ is $A_1(3m+2)$, which is divisible by $3$ because, by Chen--Ji--Jin--Shen, the three pd-rank classes of the partitions of $3m+2$ have equal size.
\end{proof}
Iterating Proposition~\ref{prop:orbit} on the base-$3$ digits of $k$ gives $F^k\equiv\prod_i F(q^{3^i})^{k_i}\pmod 3$, a bijective form of Lemma~\ref{lem:frob} for $p=3$. The same orbit argument modulo $9$ would require tracking orbits of size $3$ versus $9$ and does not give Theorem~\ref{thm:class3}; its combinatorial proof is open.

\begin{proposition}\label{prop:pdoprod}
Let $\alpha$ run over partitions into distinct parts, $\beta$ over partitions into parts $\equiv2\pmod 4$, and $\gamma$ over partitions into distinct parts $\equiv3\pmod 6$. Then $\PDO(n)$ is the number of triples $(\alpha,\beta,\gamma)$ with $|\alpha|+|\beta|+|\gamma|=n$; that is,
\[
\sum_{n\ge0}\PDO(n)\,q^n=\frac{f_4f_6^2}{f_1f_3f_{12}}=(-q;q)_\infty\,\frac1{(q^2;q^4)_\infty}\,(-q^3;q^6)_\infty .
\]
\end{proposition}
\begin{proof}
By Section~\ref{sec:prelim}, the generating function is $\prod_{n\ \mathrm{odd}}(1+q^{3n})\big/\bigl((1-q^n)(1-q^{2n})\bigr)$. Now $\prod_{n\ \mathrm{odd}}(1-q^n)^{-1}=(-q;q)_\infty$, $\prod_{n\ \mathrm{odd}}(1-q^{2n})^{-1}=1/(q^2;q^4)_\infty$, and $\prod_{n\ \mathrm{odd}}(1+q^{3n})=(-q^3;q^6)_\infty$.
\end{proof}
This is the analogue for $\mathrm{PDO}$ of the pair $(\alpha,\beta)$ of \cite{CJJS}. An explicit bijection from dotted Young diagrams with odd rows to such triples (the analogue of $\Delta$) is not constructed here, and we do not know a pd-rank for it.

\section{Conjectures}\label{sec:conj}

Computations suggest that the divisibility properties established in this paper are part of a broader, regular pattern. Specifically, our first two conjectures state that the divisibility condition derived from the initial coefficient—namely, $A_k(3) = k(k+1)(k+14)/6$ and $B_k(3) = k(k+2)(k+7)/6$—is essentially also sufficient for the progression $4n+3$.

\begin{conjecture}\label{conj:A2adic}
For every $k\ge1$, the largest power of $2$ dividing all $A_k(4n+3)$ is $2^{v_2(A_k(3))}$ if $k\not\equiv2\pmod{16}$, and exactly $2^3$ if $k\equiv2\pmod{16}$.
\end{conjecture}

\begin{conjecture}\label{conj:B2adic}
For every $k\ge1$, the largest power of $2$ dividing all $B_k(4n+3)$ is $2^{v_2(B_k(3))}$, except when $k\equiv1\pmod 8$ (where it is exactly $2^2$) or $k\equiv6\pmod 8$ (where it is exactly $2^3$, as established in Theorem~\ref{thm:pdo4exact}).
\end{conjecture}

Both statements have been verified for $k\le 96$ and the first few hundred values of $n$. Theorem~\ref{thm:mod2} and Proposition~\ref{prop:A3} establish only the lower bounds for certain residue classes of $k$. The exceptional cases are particularly striking: for $k\equiv2\pmod{16}$, $A_k(3)$ is divisible by $2^4$ or higher, yet the exact divisibility bound is precisely $2^3$.

Theorem~\ref{thm:a2_mod27} shows that for $k=2$, the progression $12n+10$ exhibits unusually deep divisibility. Among $3\le k\le27$, this phenomenon appears again only for $k=6$, where the divisibility experiences a sharp jump.

\begin{conjecture}\label{conj:A6}
For all $n\ge0$, $A_6(12n+10)\equiv0\pmod{3^5}$, and $3^5$ is the exact power of $3$ dividing all of them.
\end{conjecture}

For $k=6$, the generating function for the progression $3n+1$ can be expressed exactly as $\sum A_6(3j+1)x^j = 3\,\Theta^6\rho\,R(w)$, where $R(w)$ is a polynomial of degree $7$ in $w = x f_1^3 f_6^9 / (f_2^3 f_3^9)$ satisfying $R(w)\equiv E^6(5E-3)\pmod{81}$ with $E = 1+w$. This conjecture is equivalent to showing that certain coefficients vanish modulo $81$ in the resulting eta-quotient expansion, a fact we have not yet been able to prove. We also observe unexpectedly high divisibility in several other progressions:
\begin{conjecture}\label{conj:deep}
\begin{enumerate}
\item[\rm (a)] $A_3(12n+7) \equiv A_3(12n+11) \equiv 0 \pmod{3^4}$.
\item[\rm (b)] $A_4(6n+5) \equiv 0 \pmod{3^3}$ and $A_7(12n+11) \equiv 0 \pmod{3^3}$.
\item[\rm (c)] $A_8(18n+15) \equiv 0 \pmod{3^5}$.
\end{enumerate}
\end{conjecture}
For $k=3$, Theorem~\ref{thm:class3} yields only $9\mid A_3(3n+2)$, whereas part (a) states that divisibility jumps to $81$ along both $12n+11\subset3n+2$ and $12n+7\subset3n+1$. 

We conclude with an open question rather than a conjecture. By Corollary~\ref{cor:dict} and Proposition~\ref{prop:tri}, the divisibility $A_k(3n+2)\equiv0\pmod 3$ for $k\equiv1\pmod 3$ is equivalent to certain properties of overpartition tuples and staircase tuples; however, a combinatorial explanation remains unknown. While the pd-rank introduced in \cite{CJJS} provides such an explanation for $k=1$, simply summing the component pd-ranks weighted by powers of $\zeta=e^{2\pi i/3}$ fails for $k\ge2$.

\begin{problem}\label{prob:rank}
For $k\equiv1\pmod 3$ with $k\ge4$, find a natural statistic on $k$-tuples of designated partitions of size $3n+2$ whose values modulo $3$ divide them into three classes of equal size, preferably through an explicit bijection on Young diagrams.
\end{problem}

Solving this for $k=4$ would already be of significant interest, as $k=4$ is the first instance of the $3$-adic phenomenon in Theorem~\ref{thm:class3} where $k$ is not divisible by $3$.
\bigskip
\noindent
Taken together, our results indicate that the arithmetic of $k$-tuples of designated partitions is governed by two key components: the $3$-adic structure $F=\Theta P(u)$, which explains the progression $3n+2$ and the classes $3^b(3n+s)$ for all $k$, and a $2$-adic structure along $4n+3$ that is currently understood only partially, as reflected in Conjectures~\ref{conj:A2adic} and~\ref{conj:B2adic}. The conjectures presented here outline natural directions for future work. Ultimately, finding a crank for $k$-tuples and a combinatorial analogue of the bijection $\Delta$ for odd parts would provide the combinatorial counterparts to these two algebraic structures.
\begin{remark}
After the initial version of this preprint was posted, Hao Xiao (Universit\'e Paris-Saclay) informed us that Conjectures~6.3 and~6.4 can be proved using Radu's algorithm. The accompanying note and verification code are publicly available online.\footnote{\url{https://github.com/xiao98/pds-tuple-congruences}}
\end{remark}

\end{document}